\documentclass[a4paper, reqno, 11pt]{amsart}

\usepackage[english]{babel}
\usepackage{amsmath}
\usepackage{mathrsfs}
\usepackage{amssymb}
\usepackage{amsthm}
\usepackage{enumerate}
\usepackage{ifthen}
\usepackage{bbm}
\usepackage{color}
\usepackage{xcolor}
\usepackage{graphicx}
\provideboolean{shownotes} 
\setboolean{shownotes}{true}
\usepackage{hyperref}
\usepackage{soul}
\usepackage{geometry}
\newcommand{\margnote}[1]{
\ifthenelse{\boolean{shownotes}}%
{\marginpar{\raggedright\tiny\texttt{#1}}}%
{}%
}

\newcommand{\hole}[1]{
\ifthenelse{\boolean{shownotes}}%
{\begin{center} \fbox{ \rule {.25cm}{0cm}
\rule[-.1cm]{0cm}{.4cm} \parbox{.85\textwidth}{\begin{center}
\texttt{#1}\end{center}} \rule {.25cm}{0cm}}\end{center}}
{}
}

\newtheorem{theorem}{Theorem}[section]   

\newtheorem{proposition}[theorem]{Proposition}

\theoremstyle{definition}

\newtheorem{remark}[theorem]{Remark}

\numberwithin{equation}{section}


\begin{document}

\title[1D compressible fluids with friction]{On the no-vacuum condition for 1D compressible viscous fluids with friction}

\author[Lorenzo Pescatore]{Lorenzo Pescatore}
\address[L.Pescatore]{DISIM - Dipartimento di Ingegneria e Scienze dell'Informazione e Matematica\\ Universit\`a  degli Studi dell'Aquila \\Via Vetoio \\ 67100 L'Aquila \\ Italy}
\email[]{\href{mailto:lorenzo.pescatore@univaq.it} {lorenzo.pescatore@univaq.it}}
\begin{abstract}
For one-dimensional compressible viscous fluids, energy and BD entropy estimates prevent vacuum formation for density-dependent viscosity with exponents $0\le\alpha\le\frac12$. We show that friction $-\varepsilon u$ extends this criterion to $\alpha\le\frac23$ for both the Navier-Stokes and Navier-Stokes-Korteweg equations, assuming the strong coercivity condition $2\alpha-4\le\beta\le2\alpha-1$ in the capillary case. For the latter system, we also sharpen the no-vacuum condition $\beta\le-2$. In particular, capillarity exponents $-2<\beta<-\frac53$ are admissible under the restriction $\alpha+\beta+1\le0$, even when $\alpha>\frac23.$ Separately, for the Navier-Stokes equations with linear damping $-\varepsilon\rho u$, the active potential estimate rules out vacuum formation for every $\alpha\ge0$ and $\gamma\ge1$, provided that $\varepsilon$ is sufficiently large. Our analysis covers a broad class of power-law friction terms and, as a byproduct, yields global existence and uniqueness of strong solutions in the enlarged parameter ranges.
\end{abstract}
\maketitle
\section{Introduction}
\noindent
We are concerned with the motion of a viscous compressible fluid driven by an external force $f,$ with periodic boundary conditions
\begin{equation}\label{main system}
\begin{cases}
\partial_t \rho + \partial_x (\rho u)=0,  \\
\partial_t (\rho u)+ \partial_x (\rho u^2)=-\rho \partial_x \dfrac{\delta \mathcal{E}}{\delta \rho}+ \partial_x (\mu(\rho) \partial_x u)-\varepsilon d(\rho) u + \rho f.
\end{cases}
\end{equation}
Here $\rho=\rho(x,t)$ and $u=u(x,t)$, with $x\in\mathbb T$,
denote the density and velocity of the fluid, respectively.
The associated total energy functional is given by
\begin{equation*}
H(\rho,u)= \int_{\mathbb{T}} \bigg( \frac{1}{2} \rho u^2+ F(\rho, \partial_x \rho) \bigg) dx
\end{equation*} 
with internal energy of the fluid denoted by $$\mathcal{E}(\rho,\partial_x \rho)= \int_{\mathbb{T}} F(\rho,\partial_x \rho) dx.$$ 
We endow system \eqref{main system} with the initial conditions
\begin{equation}
(\rho,u)_{|t=0}=(\rho_0,u_0),
\end{equation}
where the initial density is nonnegative and has normalized mass.
The continuity equation yields the conservation of mass
\begin{equation*}
\int_{\mathbb T}\rho(x,t)\,dx
=\int_{\mathbb T}\rho_0(x)\,dx=1,
\end{equation*}
whereas the regularity assumptions on the initial data
will depend on the approach under consideration.
A wide class of density-dependent viscosity coefficients is considered
\begin{equation}\label{visc def}
\mu(\rho)=c_\mu \rho^\alpha, \quad c_\mu >0, \quad \alpha \ge 0,
\end{equation}
including the case of constant viscosity $\alpha=0.$
We investigate the regularizing effect of the friction term
$-\varepsilon d(\rho)u$, with $\varepsilon>0$, and its role
in preventing vacuum formation.
Besides the choices $d(\rho)=\rho$ and $d(\rho)=1$, which are widely used in the compressible fluid dynamics literature, see  \cite{BD, Bresch2, Bresch drag,Feir. Don., D.P.S.3, Feireisl damping, Lattanzio Tzavaras} and \cite{Vasseur Yu 1}, we consider the broader family of density-dependent damping coefficients
\begin{equation}\label{friction def}
d(\rho)=\rho^{\theta}, \qquad \theta \ge 0.
\end{equation}
Linear friction terms arise, for instance, in effective
descriptions of fluid motion through porous media
\cite{Whitaker}, where they model the resistance exerted
by the surrounding medium.
A density-independent resistance coefficient in the momentum
equation gives $d(\rho)=1$, as in viscous shallow-water
models with laminar bottom friction \cite{BD}, \cite{Bresch2}.
A coefficient proportional to the density, $d(\rho)=\rho$,
leads to relaxation of the velocity towards rest at a
density-independent rate.
More generally, we consider the power-law family
\eqref{friction def} as a simple model for density-dependent
resistance, investigating how the density dependence of
viscosity and friction affects the dissipative structure
of the system. Related power-law damping coefficients have also been studied for the pressureless Euler-Poisson system, where suitable density dependence can suppress singularity formation
\cite{Rozanova}.
\\
\\
The internal energy $\mathcal{E}$ is a function of the density $\rho$ and also accounts for a possible dependence on the density gradient. In particular, we consider the general framework in which the first variation of the internal energy, namely $\frac{\delta \mathcal{E}}{\delta \rho},$ is coupled with the inviscid stress tensor $\mathcal{S}$ by  
\begin{equation}\label{fund rel}
-\rho \nabla \frac{\delta \mathcal{E}}{\delta \rho}= \operatorname{div} \mathcal{S},
\end{equation} 
and depending on the expression of $\mathcal{S},$ several models are included. Below we provide a brief description of the main systems considered in our analysis.
\subsection{Navier-Stokes equations}
We consider the compressible Navier-Stokes equations in the isentropic regime. In particular, the one-dimensional inviscid stress tensor reads
\begin{equation*}
\mathcal{S}=-p(\rho)
\end{equation*}
where the quantity $p$ denotes the pressure of the fluid 
\begin{equation}\label{pressure def}
p(\rho)= c_p \rho^{\gamma}, \quad c_p >0, \quad \gamma \ge 1,
\end{equation}
and the constant $\gamma$ represents the adiabatic exponent.  The internal energy of the system \eqref{main system} can be written as
\begin{equation}
\mathcal{E}(\rho)= \int_{\mathbb{T}} h(\rho) dx, \qquad
h(\rho)=\rho \int_{1}^{\rho} \dfrac{p(z)}{z^2}dz
\end{equation}
where $h(\rho)$ is the pressure potential satisfying $\rho h''(\rho)=p'(\rho).$ 
We highlight that, in the barotropic regime given by
\eqref{pressure def}, density-dependent viscosity coefficients
have a natural physical motivation. In kinetic theory, the dependence of viscosity on temperature is determined by the microscopic interactions between molecules \cite{Chapman,Cercignani}. For hard-sphere molecules, the Chapman-Enskog viscosity law predicts that viscosity is proportional to the square root of the temperature.
Combining the ideal-gas relation with \eqref{pressure def}
then gives a viscosity proportional to $\rho^{\frac{\gamma-1}{2}}$, as recalled in \cite{Mellet}. Other molecular interaction models, including Maxwell molecules, give different relations between viscosity and temperature, providing further motivation for the power-law coefficients \eqref{visc def}. We refer to \cite{Golse} and \cite[Section 1.3.2.1]{Pescatore} for a more detailed discussion.
For completeness, our framework also includes the viscous
shallow water equations, corresponding to $\gamma=2$ and
$\alpha=1.$ 
\subsection{Navier-Stokes-Korteweg}
The Navier-Stokes-Korteweg equations describe the dynamics of a compressible fluid in which capillarity effects are not negligible. The stress tensor $\mathcal{S}$ is characterized by the additional dependence on density gradients, 
\begin{equation}
\mathcal{S}=-p(\rho)\mathbb{I}+ \mathcal{K}(\rho,\nabla \rho),
\end{equation}
encoded through the Korteweg tensor $\mathcal{K}$
\begin{equation}\label{Kort}
\mathcal{K}=\bigg( \rho \operatorname{div}(k(\rho) \nabla \rho )+ \frac{1}{2}( k(\rho)- \rho k'(\rho)) | \nabla \rho |^2 \bigg) \mathbb{I}-k(\rho) \nabla \rho \otimes \nabla \rho.
\end{equation}
The expression \eqref{Kort} is written in terms of the capillarity coefficient $k(\rho)=c_k \rho^\beta, $ with $\beta \in \mathbb{R}$ and $ c_k >0.$ In one dimension, the internal energy of system \eqref{main system} reads
\begin{equation}
\mathcal{E}(\rho, \partial_x \rho)= \int_{\mathbb{T}} \bigg( h(\rho)+ k(\rho) \dfrac{| \partial_x \rho|^2}{2} \bigg) dx.
\end{equation}
In particular, the capillary contribution represents energy
stored in spatial variations of the density.
The introduction of the Korteweg tensor is due to Korteweg
\cite{Korteweg}, while Dunn and Serrin \cite{Dunn and Serrin}
developed its thermodynamic formulation through the notion
of interstitial working.
The physical background and derivation of $\mathcal K$
are beyond the scope of this paper.
We refer to \cite[Chapter 1]{Pescatore} for an overview,
to \cite{Cahn,Cahn 0,Slemrod 0} for the diffuse-interface
description, and to \cite[Section 3.3]{Gorban} for the connection
between capillarity and kinetic theory.
\\
\\
For completeness, we also highlight a connection of the inviscid case with quantum hydrodynamics \cite{Madelung}. 
Specifically, by choosing $k(\rho)=\frac{1}{4\rho}$,
the capillarity contribution in the momentum equation becomes,
in the multidimensional case,
\begin{equation*}
\operatorname{div}\mathcal K
=\frac{1}{2}\rho\nabla
\bigg(\frac{\Delta\sqrt{\rho}}{\sqrt{\rho}}\bigg),
\end{equation*}
which corresponds to the Bohm quantum force.
This term appears in quantum hydrodynamic models used
in the analysis of semiconductors and Bose-Einstein condensates; see, for instance \cite{Landau}, \cite{Antonelli Marcati} and \cite{Ant4}.
\subsection{No vacuum analysis $\varepsilon=0$}
Here we provide a brief description of the main elements used in order to prevent the formation of vacuum states of the density for one-dimensional viscous fluids. These techniques naturally exhibit some regimes of the constitutive parameters which we aim to improve by considering the influence of a friction term.
\\
\\
\underline{\textit{BD Entropy:}}
A key feature in the study of viscous compressible models is the use of a Bresch-Desjardins entropy method. This analysis relies on the effective velocity $V=u+Q,$ see for instance \cite{BD},\cite{Bresch2},\cite{Const}, \cite{Haspot 2}, \cite{Mellet}, where $Q$ satisfies the transport identity
$$Q:= \mu(\rho) \dfrac{\partial_x \rho}{\rho^2}, \quad \partial_t Q+ u \partial_x Q=-\dfrac{1}{\rho} \partial_x(\mu(\rho)\partial_x u).$$
This change of variables allows one to rewrite system \eqref{main system} as a system of balance  laws in the new variables $(\rho, \rho V), $ see \eqref{balance laws} and \eqref{balance laws NSK}.  The BD entropy estimate then represents an energy estimate for the new system, involving the functional 
\begin{equation}\label{BD Functional}
H^{BD}(\rho,u)= \int_{\mathbb{T}} \bigg( \frac{1}{2} \rho |u+Q|^2+ F(\rho, \nabla \rho) \bigg) dx.
\end{equation}
See also \cite{Const2}, where a hierarchy of higher-order
entropies is developed for one-dimensional compressible
models with local and nonlocal dissipation operators.
\\
In the absence of capillarity,  namely $c_k=0,$ the internal energy density $F$ is a function of the density only. Therefore the energy $H(\rho,u)$ provides no information on the spatial regularity of density gradients. Conversely, the energy estimate and the bound on the BD entropy functional \eqref{BD Functional} control $\sqrt{\rho}Q$, and hence $\partial_x\rho^{\alpha-\frac12}$ in
$L^\infty(0,T;L^2(\mathbb T)),$ see Proposition
\ref{Prop NS BD}.
In one dimension, this density-gradient control allows us
to bound $\|\rho^{-1}\|_{L^\infty}$ in terms of the initial
data through Sobolev embedding.
This yields a no-vacuum condition provided that 
\begin{equation*}\label{NV NS} \tag*{(NV)$_{\mathrm{NS}}$}
0 \le \alpha \le \dfrac{1}{2}.
\end{equation*}
In the case of capillary fluids, the derivation of the entropy estimate is more involved due to the interaction of the Korteweg tensor with the transport term $Q$, see Proposition \ref{Prop NSK BD}. In particular, the tensor $\mathcal{K}$ is accounted for general capillarity exponents $\beta \in \mathbb{R}$ and several nonlinear terms appear in the analysis.  On the one hand, the validity of the BD entropy method for capillary fluids is obtained in the regime of the parameters $\alpha$ and $\beta$ dictated by the \textit{strong coercivity condition} \cite{Le Floch}
\begin{equation}\label{scc}
2\alpha-4 \le \beta \le 2\alpha-1,
\end{equation}
which characterizes the non-negativity of an entropy dissipation functional $\mathcal{F}_{k}[\rho],$ see Proposition \ref{Prop NSK BD}. On the other hand, already the energy estimate yields $\partial_x \rho^{\frac{\beta}{2}+1} \in L^\infty(0,T; L^2(\mathbb{T})).$ Thus, under \eqref{scc}, the energy and BD entropy estimates
yield the following no-vacuum criterion, based on the
no-cavitation mechanism studied in \cite{Le Floch}:
\begin{equation*}\label{NV NSK} \tag*{(NV)$_{\mathrm{NSK}}$}
0 \le \alpha \le \dfrac{1}{2}, \quad \text{or} \quad \beta \le -2.
\end{equation*}
At the endpoints $\alpha=\frac12$ and $\beta=-2$, the preceding
density-gradient estimates are understood in logarithmic form,
involving $\partial_x\log\rho.$
\\
\\
\underline{\textit{Effective pressure and active potential:}}
For the one-dimensional Navier-Stokes equations, further information on the initial data allows one to obtain lower density bounds beyond the range covered by the BD entropy estimate. Suppose that
\begin{equation}\label{reg initial datum density}
0<\underline{\rho}\le\rho_0\le\overline{\rho}.
\end{equation}
For the undamped Cauchy problem on $\mathbb R$, Burtea and Haspot~\cite{BurteaHaspot} introduce the effective pressure, which in our notation takes the form
\begin{equation*}
Y=\frac{\partial_xV}{\rho}+\int_1^\rho\frac{p'(s)}{s\mu(s)}\,ds.
\end{equation*}
For $\alpha>\frac12$ and $\gamma\ge\max\{1,\alpha\}$, Burtea and Haspot use an upper bound on $Y$ to control $\rho^{-1}$ on every finite time interval. Their argument requires a one-sided bound on the initial effective velocity derivative, which on the torus reads $\partial_xV_0\le K_0\rho_0$ for some $K_0\ge0$. On the one hand, this is an extra assumption at the level of the energy and BD entropy. On the other hand, for the Sobolev data considered in Section~\ref{Sec global wp}, it holds automatically since $\partial_xV_0$ is bounded and $\rho_0$ is bounded away from zero.
\\
\\
Another way to obtain a pointwise bound on the density is to use the active potential introduced in \cite{Const},
\begin{equation}
w=-p(\rho)+\mu(\rho)\partial_xu.
\end{equation}
This is the one-dimensional effective viscous flux \cite{Feir Nov Petz,Lions}, whose spatial derivative appears in the momentum equation. For a spatially homogeneous force $f=f(t)$, the reaction term in the equation for $w$ preserves the sign $w\le0$ whenever $w_0\le0$ and $\gamma\le\alpha+1$. If in addition $\gamma\ge\alpha$, the continuity equation then gives a positive lower bound for the density on every finite time interval. We refer to \cite[Section~6]{Const} for the argument.
\\
\\
For Korteweg fluids, including the capillary stress in the definition of the active potential,
\begin{equation*}
w=-p(\rho)+\mu(\rho)\partial_xu+\mathcal K,
\end{equation*}
leads to an evolution equation involving fourth-order spatial derivatives of the density. The preceding maximum-principle argument therefore does not extend directly to this setting.
\subsection{Friction regularization $\varepsilon >0$}
The first main goal of the present paper is to improve the
preceding energy/entropy analysis by exploiting friction.
Its interaction with the effective velocity $V=u+Q$
generates an additional density functional in the BD entropy
balance.
More precisely, we introduce the friction potential
$I(\rho;\varepsilon)$ through
\begin{equation}\label{friction potential}
I_{\rho\rho}(\rho;\varepsilon)
=\varepsilon\frac{\mu(\rho)d(\rho)}{\rho^3},
\qquad
I(1;\varepsilon)=I_\rho(1;\varepsilon)=0.
\end{equation}
For smooth positive solutions, the continuity equation yields
\begin{equation*}
\varepsilon\int_{\mathbb T}d(\rho)uV\,dx
=\varepsilon\int_{\mathbb T}d(\rho)u^2\,dx
+\frac{d}{dt}\int_{\mathbb T}I(\rho;\varepsilon)\,dx.
\end{equation*}
Interestingly, a density functional with a comparable structure arises
in the $\kappa$-entropy analysis of density-independent friction
by Bresch, Vasseur and Yu \cite[Section~2.2.2]{Bresch}. The potential $I$ is convex and nonnegative. For the energy and BD entropy estimates, we assume
$f\in L^2(0,T;L^\infty(\mathbb T))$ and
\begin{equation}\label{reg initial datum}
H(\rho_0,u_0)<\infty,
\qquad
H^{BD}(\rho_0,u_0)
+\int_{\mathbb T}I(\rho_0;\varepsilon)\,dx<\infty.
\end{equation}
For the power-law coefficients
\eqref{visc def} and \eqref{friction def},
$I$ controls negative powers of the density whenever
$\alpha+\theta<1$; see Remark
\ref{remark NS vacuum sec 2}.
Combining this control with the spatial regularity supplied
by the energy and BD entropy estimates yields the following
improved no-vacuum conditions for $d(\rho)=1$:
\begin{equation*}\label{epsilon NV NS}
\tag*{(NV)$_{\mathrm{NS}_{\varepsilon}}$}
0\le\alpha\le\frac23,
\end{equation*}
\begin{equation*}\label{epsilon NV NSK}
\tag*{(NV)$_{\mathrm{NSK}_{\varepsilon}}$}
2\alpha-4\le\beta\le2\alpha-1,
\qquad
0\le\alpha\le
\max\left\{\frac23,-\beta-1\right\}.
\end{equation*}
For the general case $\theta\ge0$, we refer to the no vacuum conditions in Propositions \ref{Prop density NS} and \ref{Prop density NSK}.
\\
Damped systems are also used as intermediate approximations
in the construction of global weak solutions
\cite{Bresch,Vasseur Yu 0,Vasseur Yu 1}.
In particular, the density-independent drag in
\cite{Vasseur Yu 1} provides logarithmic control of the
density, excluding vacuum sets of positive measure.
In future work, we aim to exploit the improved conditions
obtained through the BD entropy method to construct
approximations with strictly positive density for a broader
class of viscosity and capillarity coefficients.
\\
\\
We next consider the pointwise estimates for the Navier-Stokes equations. In Proposition~\ref{Prop effective pressure NS} we extend the effective-pressure argument of Burtea and Haspot to a class of power-law friction terms. For $\gamma\ge\alpha$, this gives a lower density bound on every finite time interval under the one-sided condition $\partial_xV_0\le K_0\rho_0$. For linear damping $-\varepsilon\rho u$ and a spatially homogeneous force, the pressure dissipation in the BD entropy yields the stronger estimate
\begin{equation*}
\rho(x,t)\ge c>0
\qquad\text{for all }x\in\mathbb T,\quad 0\le t<T,
\end{equation*}
where $c$ is independent of $T$. This uniform bound also holds in the undamped case under the same initial one-sided condition.
Linear damping also contributes the term $-\varepsilon(w+p(\rho))$ to the active-potential equation. For every $\varepsilon>0$ and $\gamma\le\alpha+1$, we obtain an explicit time $T_\varepsilon$ such that, on any interval $[0,T)$ of smooth existence with $T>T_\varepsilon$,
\begin{equation*}
w(x,t)\le0
\qquad\text{for all }x\in\mathbb T,\quad T_\varepsilon\le t<T,
\end{equation*}
regardless of the initial sign of $w$. The role of damping is particularly clear at $\gamma=\alpha+1$, where the reaction at $w=0$ reduces to $-\varepsilon p(\rho)$. In the strict range $\gamma<\alpha+1$, the pressure term can itself force $w$ to become nonpositive even without friction.
\\
\\
In the range $\alpha\le\gamma\le\alpha+1$, the two pointwise arguments can be used together: the effective pressure controls the density up to $T_\varepsilon$, and the sign of $w$ gives a lower bound uniform in time thereafter. If $w_0\le0$, the latter argument applies from the initial time and the one-sided condition on $\partial_xV_0$ is unnecessary. We refer to Proposition~\ref{prop:ap-sign-entry} for both conclusions.
\\
\\
Finally, for linear damping and $f(x,t)=f(t),$ sufficiently large $\varepsilon$ gives a quantitative bound on both density extrema for every $\alpha\ge0$ and $\gamma\ge1$. More precisely, writing $m_0=\min_{\mathbb T}\rho_0$ and $M_0=\max_{\mathbb T}\rho_0$, we prove that
\begin{equation*}
m_0-\frac{C}{\varepsilon}\le\rho(x,t)\le M_0+\frac{C}{\varepsilon},
\end{equation*}
where $C$ is independent of time and of $\varepsilon$. Given $0<L<m_0\le M_0<R$, sufficiently large damping keeps the density in $[L,R]$ from the initial time.  As discussed in Remark~\ref{remark high-friction}, the bound on the density extrema can be viewed as a finite-friction counterpart of the maximum principle for the limiting diffusion equation. For the damped Euler equations, the high-friction limit was justified in \cite{Lattanzio Tzavaras}, and its connection with gradient flows is discussed in \cite{Lattanzio Tzavaras CPDE}. The unforced Navier-Stokes equations have the same formal limit, since viscosity enters at higher order under this scaling, as observed in \cite{CianfaraniLattanzio} for capillary fluids.
\\
\\
The resulting global existence and uniqueness theorems for the Navier-Stokes and Navier-Stokes-Korteweg equations with friction are stated in Theorem \ref{thm:NS-global} and Theorem \ref{thm:NSK-global}, respectively.
\subsection*{State of the art}
We briefly recall some existence and uniqueness results for the models considered above in the absence of friction.
For the compressible Navier-Stokes equations with constant viscosity, we refer to the pioneering work of Lions and of Feireisl and collaborators on finite-energy weak solutions \cite{Lions,Feir Nov Petz}. For the one-dimensional Cauchy problem with constant viscosity, Hoff~\cite{Hoff1998} constructed global weak solutions for large, possibly discontinuous data. He also showed that if the initial density is bounded above and away from zero, it remains so on every finite time interval.
For density-dependent viscosities satisfying the BD structure, global weak solutions allowing vacuum were constructed in \cite{Bresch,Li Xin,Vasseur Yu 0}.
In one space dimension, Mellet and Vasseur \cite{Mellet} constructed global strong solutions for large initial data with density bounded away from zero, allowing the viscosity to degenerate at vacuum.
Haspot \cite{Haspot 2} proved global existence of strong solutions for power-law viscosities, including the viscous shallow-water case, using estimates for the effective velocity and a maximum principle for the density. For results on finite-time blow-up in one-dimensional compressible Navier-Stokes flows, we refer to \cite{LiLiXin,CaoLiYu}.
In the presence of vacuum, local strong well-posedness for three-dimensional fluids with nondegenerate viscosities was established in \cite{Cho 2004}. See also \cite{Don-Pes} for local solvability and blow-up criteria for related models with external forces.
In two dimensions, Danchin and Mucha \cite{DanchinMucha} constructed suitable global weak solutions with bounded density and large bulk viscosity, and proved uniqueness for the linear pressure law. For a recent uniqueness result covering more general pressure laws, we refer to \cite{PescatoreViolini}.
\\
\\
For capillary fluids, Lacroix-Violet and Vasseur \cite{Lacroix} and Antonelli and Spirito \cite{Spirito1} proved global existence of weak solutions to the quantum Navier-Stokes equations in two and three dimensions, allowing for vacuum.
For a three-dimensional Korteweg model with linear viscosity and constant capillarity, we refer to \cite{Spirito}.
In one dimension, Germain and LeFloch \cite{Le Floch} developed a finite-energy framework covering several classes of viscosity and capillarity coefficients, including regimes with vacuum. Antonelli, Bresch and Spirito \cite{Ant} subsequently enlarged the power-law regimes admitting global weak solutions with vacuum, provided that $2\alpha-3 \le \beta < 2\alpha-1,$ with $\alpha>\frac{1}{2}$ and $\beta >-2.$ For quantum capillarity, related weak existence results were obtained in \cite{D.P.S.2}.
\\
Concerning strong solutions, we refer to \cite{Chen} for
global classical solvability for large data away from vacuum,
under restrictions on the viscosity and capillarity exponents,
and to \cite{Charve,Burtea1,Burtea2} for global solvability
and vanishing-capillarity results under specific relations
between the two coefficients.
Local well-posedness for every $\alpha\ge0$ and
$\beta\in\mathbb R$, together with global strong solvability
under \eqref{scc} and \ref{NV NSK}, was established
in \cite{Pesc}.
Related strong-solution results for quantum capillarity
were obtained in \cite{D.P.S.}.
We highlight that the results \cite{Pesc,D.P.S.,D.P.S.2} account for both the deterministic and stochastically forced cases.
\subsection*{Organization of the paper}
The outline of the paper is the following.  Section 2 is devoted to the BD entropy analysis for both the Navier-Stokes and Navier-Stokes-Korteweg equations. In Section \ref{sec:active potential}, we first use the effective pressure to obtain lower density bounds for a class of power-law friction terms. We then study the active potential for the Navier-Stokes equations with linear damping. Finally, Section \ref{Sec global wp} is devoted to the proof of two global well-posedness results, Theorem \ref{thm:NS-global} and Theorem \ref{thm:NSK-global} in the enlarged range of constitutive parameters.
\subsection*{Acknowledgements}
The author gratefully acknowledges the partial support by the Gruppo Nazionale per l’Analisi Matematica, la Probabilità e le loro Applicazioni (GNAMPA) of the Istituto Nazionale di Alta Matematica (INdAM).
\section{A priori estimates}\label{Sec2}
\noindent
We derive the energy and BD entropy estimates for both systems and use them to obtain upper and lower bounds on the density.
\subsection*{Notation and conventions.} Throughout the paper, we assume the following conventions.  The Lebesgue space is denoted by $L^p(\mathbb{T})$ with related norm $\| \cdot \|_{L^p}.$ Similarly $W^{k,p}(\mathbb{T})$ denotes the Sobolev space of $L^p(\mathbb{T})$ functions with $k$ distributional derivatives in $L^p(\mathbb{T})$ and $H^k(\mathbb{T})$ corresponds to the case $p=2.$ For a Banach space $X$, we denote by $C([0,T];X)$
the space of continuous $X$-valued functions and by
$L^p(0,T;X)$ the usual Bochner space.  For simplicity, we assume $|\mathbb{T}|=1$ and that the total mass is normalized to one. 
\\
\\
In what follows, we denote by $C=C(T)$ a generic constant which depends on $T,f$, energy/entropy norms of the initial conditions $(\rho_0,u_0)$, the initial friction potential $\int_{\mathbb T}I(\rho_0;\varepsilon)\,dx$, and possibly on $\varepsilon$, see Proposition \ref{Prop NS BD} and Proposition
\ref{Prop NSK BD}. The constants may vary from line to line and do not depend on a lower bound assumed for the density along the solution. Throughout this section, we derive the estimates on intervals where the solution is smooth and the density is strictly positive.
\subsection{Navier-Stokes equations}
We consider the isentropic Navier-Stokes equations
\begin{equation}\label{NS damping}
\begin{cases}
\partial_t\rho+\partial_x(\rho u)=0,\\
\partial_t(\rho u)+\partial_x(\rho u^2)
=-\partial_xp(\rho)+\partial_x(\mu(\rho)\partial_xu)
-\varepsilon d(\rho)u+\rho f,
\end{cases}
\end{equation}
for which the internal energy $\mathcal E$ reads
\begin{equation*}
\mathcal E(\rho)=\int_{\mathbb T}h(\rho)\,dx,
\qquad
h(\rho)=\rho\int_1^\rho\frac{p(z)}{z^2}\,dz.
\end{equation*}
We start our analysis with the following standard energy estimate.
\begin{proposition}\label{Prop NS Energy}
Let $(\rho,u)$ be a smooth solution of \eqref{NS damping} with finite
initial energy. Then the following inequality holds:
\begin{equation}\label{NS energy ineq}
\begin{split}
\sup_{t\in[0,T]}H(\rho,u)(t)
&+\int_0^T\int_{\mathbb T}\mu(\rho)|\partial_xu|^2\,dxdt+\varepsilon\int_0^T\int_{\mathbb T}d(\rho)u^2\,dxdt
\le C(T).
\end{split}
\end{equation}
\end{proposition}

\begin{proof}
The proof easily follows by testing the momentum equation
$\eqref{NS damping}_2$ with $u$, integrating by parts and using the
continuity equation. For completeness, we observe that the external
force $f$ gives the following supply rate in the energy balance:
\begin{equation*}
\begin{split}
\int_{\mathbb T}\rho fu\,dx
&\le\frac12\int_{\mathbb T}\rho f^2\,dx
+\frac12\int_{\mathbb T}\rho u^2\,dx\le\frac12\|f\|_{L_x^\infty}^2
+\frac12\int_{\mathbb T}\rho u^2\,dx.
\end{split}
\end{equation*}
Therefore \eqref{NS energy ineq} follows as a standard application of
the Gronwall lemma.
\end{proof}

\noindent
The next proposition is devoted to the BD entropy estimate. To begin with,
we introduce the effective velocity $V$:
\begin{equation*}
V=u+Q,
\qquad
Q=\mu(\rho)\frac{\partial_x\rho}{\rho^2}.
\end{equation*}
The quantity $Q$ solves the following inhomogeneous transport equation:
\begin{equation*}
\partial_tQ+u\partial_xQ
=-\frac1\rho\partial_x(\mu(\rho)\partial_xu).
\end{equation*}
Thus we may rewrite system \eqref{NS damping} in terms of the new unknowns
$(\rho,\rho V)$:
\begin{equation}\label{balance laws}
\begin{cases}
\partial_t\rho+\partial_x(\rho u)=0,\\
\partial_t(\rho V)+\partial_x(\rho uV)
=-\partial_xp(\rho)-\varepsilon d(\rho)u+\rho f
\end{cases}
\end{equation}
and we now derive the BD entropy estimate, including the contribution of the friction potential introduced in \eqref{friction potential}. Note that the quantity $\frac{\partial_x\rho^m}{m}$ is understood as $\partial_x\log\rho$ in the degenerate case $m=0.$
\begin{proposition}\label{Prop NS BD}
Let $(\rho,u)$ be a smooth solution of \eqref{NS damping} with initial
datum satisfying \eqref{reg initial datum}. Then the following inequality holds:
\begin{equation}\label{NS BD ineq}
\begin{split}
\sup_{t\in[0,T]}
\left(
H^{BD}(\rho,u)(t)
+\int_{\mathbb T}I(\rho;\varepsilon)\,dx
\right) & + \frac{4\gamma c_pc_\mu}{(\gamma+\alpha-1)^2}
\int_0^T\int_{\mathbb T}
\left|\partial_x\rho^{\frac{\alpha+\gamma-1}{2}}\right|^2
\,dxdt \\ & +\varepsilon\int_0^T\int_{\mathbb T}d(\rho)u^2\,dxdt
\le C(T),
\end{split}
\end{equation}
with the friction potential defined in \eqref{friction potential}.
\end{proposition}

\begin{proof}
We test $\eqref{balance laws}_2$ with $V=u+Q$, and after integrating
by parts we have
\begin{equation*}
\int_{\mathbb T}
\left[\partial_t(\rho V)+\partial_x(\rho uV)\right](u+Q)\,dx
=\frac{d}{dt}\int_{\mathbb T}\frac12\rho|u+Q|^2\,dx.
\end{equation*}
The pressure term gives rise to the following entropy dissipation
functional:
\begin{equation*}
\begin{split}
\int_{\mathbb T}\partial_xp(\rho)Q\,dx
&=\gamma c_pc_\mu
\int_{\mathbb T}\rho^{\gamma+\alpha-3}|\partial_x\rho|^2\,dx=\frac{4\gamma c_pc_\mu}{(\gamma+\alpha-1)^2}
\left\|\partial_x\rho^{\frac{\gamma+\alpha-1}{2}}\right\|_{L^2}^2.
\end{split}
\end{equation*}
The contribution due to the friction term can be handled by using the
continuity equation:
\begin{equation*}
\begin{split}
\varepsilon\int_{\mathbb T}d(\rho)uQ\,dx
&=\varepsilon c_\mu
\int_{\mathbb T}\rho^{\alpha+\theta-2}u\partial_x\rho\,dx=\int_{\mathbb T}
\rho I_{\rho\rho}(\rho;\varepsilon)u\partial_x\rho\,dx=\frac{d}{dt}\int_{\mathbb T}I(\rho;\varepsilon)\,dx.
\end{split}
\end{equation*}
Since
\[
I_{\rho\rho}(\rho;\varepsilon)
=\varepsilon c_\mu\rho^{\alpha+\theta-3}\ge0
\]
and $I(1;\varepsilon)=I_\rho(1;\varepsilon)=0$, the function $I$ is
convex and nonnegative. Therefore
\begin{equation}\label{lower bound I}
\int_{\mathbb T}I(\rho;\varepsilon)\,dx\ge0.
\end{equation}
For completeness, we also highlight that
\begin{equation*}
\begin{split}
\int_{\mathbb T}\rho f(u+Q)\,dx
&\le\frac12\int_{\mathbb T}\rho f^2\,dx
+\frac12\int_{\mathbb T}\rho|u+Q|^2\,dx \le\frac12\|f\|_{L_x^\infty}^2
+\frac12\int_{\mathbb T}\rho|u+Q|^2\,dx,
\end{split}
\end{equation*}
while the remaining terms involving only the velocity $u$ are treated
as in the energy inequality, Proposition \ref{Prop NS Energy}.
The proof of \eqref{NS BD ineq} then follows by using the Gronwall
lemma applied to the quantity
\begin{equation*}
Y(t)=H^{BD}(\rho,u)(t)
+\int_{\mathbb T}I(\rho;\varepsilon)\,dx
\end{equation*}
and recalling \eqref{lower bound I}.
In particular, we highlight that
\begin{equation}\label{NS V estimate}
\sup_{t\in[0,T]}
\int_{\mathbb T}\frac12\rho|V|^2\,dx\le C(T)
\end{equation}
which combined with Proposition \ref{Prop NS Energy} gives
\begin{equation}\label{NS Q estimate}
\begin{split}
\sup_{t\in[0,T]}\int_{\mathbb T}\rho Q^2\,dx \le C(T).
\end{split}
\end{equation}
\end{proof}
\begin{remark}\label{remark NS vacuum sec 2}
It is important to mention that by definition of the damping potential $I(\rho; \varepsilon)$ and Proposition \ref{Prop NS BD}, we have the following $\varepsilon$-dependent bound
\begin{equation*}
\varepsilon\sup_{t\in[0,T]}
\int_{\mathbb T}\rho^{\alpha+\theta-1}\,dx\le C(T),
\end{equation*}
which gives additional integrability properties of negative powers
of the density whenever $\alpha+\theta<1.$
\end{remark}
\noindent
We conclude this short subsection with the derivation of suitable
bounds for the density, including a no-vacuum condition.
\begin{proposition}\label{Prop density NS}
Let $(\rho,u)$ be a smooth solution of \eqref{NS damping} with initial
datum satisfying \eqref{reg initial datum}.
Then the following upper and lower bounds for the density hold:
\begin{itemize}
\item[1)] For any $\alpha$ and $\theta$,
\begin{equation}\label{upper rho NS}
\rho\in L^\infty((0,T)\times\mathbb T).
\end{equation}
\item[2)] For any $\alpha$ and $\theta$ such that
$0\le\alpha\le\max\{\frac12,\frac{2-\theta}{3}\}$,
\begin{equation}\label{no vacuum NS}
\frac1\rho\in L^\infty((0,T)\times\mathbb T).
\end{equation}
\end{itemize}
\end{proposition}

\begin{remark}\label{remark vacuum NS}
The estimate \eqref{no vacuum NS} suggests that the friction
regularization holds whenever $\theta\in[0,\frac12)$.
In particular, for $\theta=0$ the no-vacuum condition reduces to
$0\le\alpha\le\frac23$, which is exactly \ref{epsilon NV NS}.
\end{remark}

\begin{proof}
To begin with, we observe that by virtue of the conservation of mass
and the mean value theorem, for each $t$ there exists a point
$x_0\in\mathbb T$ such that $\rho(x_0,t)=1$.
For $\alpha>0$, the fundamental theorem of calculus, H\"older
inequality, Proposition \ref{Prop NS Energy} and Proposition \ref{Prop NS BD} yield
\begin{equation}\label{tfc NS}
\begin{split}
\rho^\alpha(x,t)
&=1+\int_{x_0}^x\partial_y\rho^\alpha(y,t)\,dy \le1+\frac{\alpha}{c_\mu}
\int_{\mathbb T}\sqrt\rho\,|\sqrt\rho Q|\,dx
\le C.
\end{split}
\end{equation}
Therefore \eqref{upper rho NS} follows by taking the supremum in
space and time. In the case $\alpha=0$, we instead have
\begin{equation*}
|\log\rho(x,t)|
\le\frac1{c_\mu}
\int_{\mathbb T}\sqrt\rho\,|\sqrt\rho Q|\,dx
\le C.
\end{equation*}
Concerning the no-vacuum bound \eqref{no vacuum NS}, first we notice
that by virtue of \eqref{NS Q estimate}, Sobolev embedding and the
identity $\rho(x_0,t)=1$, we have
$\rho^{\alpha-\frac12}\in L^\infty((0,T)\times\mathbb T)$,
which yields \eqref{no vacuum NS} for $\alpha<\frac12$.
At $\alpha=\frac12$, the estimate \eqref{NS Q estimate} controls
$c_\mu\partial_x\log\rho$, and hence
\[
|\log\rho(x,t)|
\le\int_{\mathbb T}|\partial_x\log\rho|\,dx\le C.
\]
Additionally, we observe that
\begin{equation*}
\begin{split}
|\log\rho(x,t)|
&\le\int_{\mathbb T}|\partial_x\log\rho|\,dx=\frac1{|\alpha-\frac12|}
\int_{\mathbb T}
\rho^{\frac{2-3\alpha-\theta}{2}}
\rho^{\frac{\alpha+\theta-1}{2}}
\left|\partial_x\rho^{\alpha-\frac12}\right|\,dx
\le C,
\end{split}
\end{equation*}
provided that $\alpha>\frac12$ and $2-3\alpha-\theta\ge0$.
These conditions imply $\alpha+\theta<1.$ 
\end{proof}
\subsection{Navier-Stokes-Korteweg equations}
We complete this section with the analysis of the compressible
Navier-Stokes-Korteweg equations. To begin with, we recall the
system
\begin{equation}\label{NSK damping}
\begin{cases}
\partial_t\rho+\partial_x(\rho u)=0,\\
\partial_t(\rho u)+\partial_x(\rho u^2)
=-\partial_xp(\rho)+\partial_x(\mu(\rho)\partial_xu)
+\partial_x\mathcal K-\varepsilon d(\rho)u+\rho f,
\end{cases}
\end{equation}
with Korteweg tensor expressed by
\begin{equation*}
\mathcal K(\rho,\partial_x\rho)
=\rho\partial_x(k(\rho)\partial_x\rho)
-\frac12\big(k(\rho)+\rho k'(\rho)\big)|\partial_x\rho|^2,
\end{equation*}
for $k(\rho)$ being the capillarity coefficient satisfying the
following power-law relation:
\begin{equation*}
k(\rho)=c_k\rho^\beta,
\qquad \beta\in\mathbb R,
\qquad c_k>0.
\end{equation*}
It is worth mentioning that the third-order derivative operator
arising in the momentum equation can be rewritten as
\begin{equation*}
\partial_x\mathcal K(\rho,\partial_x\rho)
=\rho\partial_x\left(
\partial_x(k(\rho)\partial_x\rho)
-\frac{k'(\rho)}2|\partial_x\rho|^2
\right),
\end{equation*}
while the internal energy is defined as
\begin{equation}
\mathcal E(\rho,\partial_x\rho)
=\int_{\mathbb T}
\left(h(\rho)+k(\rho)\frac{|\partial_x\rho|^2}{2}\right)\,dx.
\end{equation}
By an argument similar to that used in the zero-capillarity case,
we have the following standard energy inequality.

\begin{proposition}\label{Prop NSK Energy}
Let $(\rho,u)$ be a smooth solution of \eqref{NSK damping} with finite
initial energy. Then the following inequality holds:
\begin{equation}\label{NSK energy ineq}
\begin{split}
\sup_{t\in[0,T]}H(\rho,u)(t)
&+\int_0^T\int_{\mathbb T}\mu(\rho)|\partial_xu|^2\,dxdt +\varepsilon\int_0^T\int_{\mathbb T}d(\rho)u^2\,dxdt
\le C(T).
\end{split}
\end{equation}
\end{proposition}

\begin{proof}
We only consider the new terms arising due to the capillarity term
in the momentum equation. Specifically, by multiplying
$\eqref{NSK damping}_2$ by $u$ and integrating by parts we have
\begin{equation*}
\begin{split}
\int_{\mathbb T}u\,\partial_x\mathcal K\,dx
&=\int_{\mathbb T}
\partial_x\left(
\partial_x(k(\rho)\partial_x\rho)
-\frac{k'(\rho)}2|\partial_x\rho|^2
\right)\rho u\,dx=-\frac{d}{dt}\int_{\mathbb T}
k(\rho)\frac{|\partial_x\rho|^2}{2}\,dx.
\end{split}
\end{equation*}
\end{proof}

\noindent
We emphasize that for Korteweg fluids, the energy bounds
involve the gradient of the density
$\partial_x\rho^{\frac\beta2+1}\in L^\infty_tL^2_x$ when
$\beta\neq-2$, and
$\partial_x\log\rho\in L^\infty_tL^2_x$ when $\beta=-2$.
Hence the vacuum analysis also relies on the capillarity exponent
$\beta$. Contrary to the zero-capillarity case, here the BD entropy
method holds under suitable restrictions on $\alpha$ and $\beta.$ Specifically, we state and prove the following proposition.  To aid the reader, we also refer to \cite[Proposition 2.3]{Le Floch} and \cite[Proposition 2.3]{Ant} for a similar analysis without friction.
\begin{proposition}\label{Prop NSK BD}
Let $(\rho,u)$ be a smooth solution of \eqref{NSK damping} with initial datum satisfying \eqref{reg initial datum}. Assume in addition that \eqref{scc} is satisfied.
Then the following inequality holds:
\begin{equation}\label{NSK BD ineq}
\begin{split}
\sup_{t\in[0,T]}
\left(
H^{BD}(\rho,u)(t)
+\int_{\mathbb T}I(\rho;\varepsilon)\,dx
\right) & + \frac{4\gamma c_pc_\mu}{(\gamma+\alpha-1)^2}
\int_0^T\int_{\mathbb T}
\left|\partial_x\rho^{\frac{\alpha+\gamma-1}{2}}\right|^2
\,dxdt \\ & +\int_0^T\mathcal F_k[\rho]\,dt
+\varepsilon\int_0^T\int_{\mathbb T}d(\rho)u^2\,dxdt
\le C(T),
\end{split}
\end{equation}
where for $\alpha+\beta+1\neq0$, the entropy dissipation functional is
\begin{equation}
\mathcal F_k[\rho]
=\frac{4c_\mu c_k}{(\alpha+\beta+1)^2}
\int_{\mathbb T}
\left|\partial_{xx}\rho^{\frac{\alpha+\beta+1}{2}}\right|^2\,dx
+c(\alpha,\beta)\int_{\mathbb T}
\left|\partial_x\rho^{\frac{\alpha+\beta+1}{4}}\right|^4\,dx,
\end{equation}
with
\begin{equation*}
c(\alpha,\beta)
=\frac{64c_\mu c_k}{(\alpha+\beta+1)^2}
\left[
\frac{(\alpha-\beta-1)(1-\alpha)}{(\alpha+\beta+1)^2}
-\frac{\beta}{3(\alpha+\beta+1)}
\right],
\end{equation*}
while if $\alpha+\beta+1=0$, it is given by
\begin{equation}
\mathcal F_k[\rho]
=c_\mu c_k\int_{\mathbb T}|\partial_{xx}\log\rho|^2\,dx
+\frac{c_\mu c_k\alpha(1-\alpha)}2
\int_{\mathbb T}|\partial_x\log\rho|^4\,dx.
\end{equation}
\end{proposition}

\begin{proof}
First we rewrite system \eqref{NSK damping} in the balance form
\begin{equation}\label{balance laws NSK}
\begin{cases}
\partial_t\rho+\partial_x(\rho u)=0,\\
\partial_t(\rho V)+\partial_x(\rho uV)
=-\partial_xp(\rho)+\partial_x\mathcal K
-\varepsilon d(\rho)u+\rho f.
\end{cases}
\end{equation}
Then we test $\eqref{balance laws NSK}_2$ with $V=u+Q$ and we focus
on the analysis of the key quantity
\begin{equation*}
\int_{\mathbb T}\partial_x\mathcal K\,Q\,dx
=c_\mu\int_{\mathbb T}
\partial_x\left(
\partial_x(k(\rho)\partial_x\rho)
-\frac{k'(\rho)}2|\partial_x\rho|^2
\right)\rho^{\alpha-1}\partial_x\rho\,dx.
\end{equation*}
After a tedious rearrangement of density derivatives, we integrate
by parts to obtain
\begin{equation*}
\begin{split}
\int_{\mathbb T}\partial_x\mathcal K\,Q\,dx
=& -c_\mu c_k\int_{\mathbb T}\bigg[
\rho^{\alpha+\beta-1}|\partial_{xx}\rho|^2+\left(\alpha+\frac\beta2-1\right)
\rho^{\alpha+\beta-2}|\partial_x\rho|^2\partial_{xx}\rho\bigg] dx \\ & -c_{\mu} c_k \int_{\mathbb{T}} \bigg[\frac{\beta(\alpha-1)}2
\rho^{\alpha+\beta-3}|\partial_x\rho|^4
\bigg] \,dx =-\mathcal F_k[\rho].
\end{split}
\end{equation*}
On the one hand, the constant $c(\alpha,\beta)$ is allowed to take
negative values. On the other hand, it has already been pointed out
in several contributions that $\mathcal F_k[\rho]\ge0$ if and only if \eqref{scc} is satisfied. Specifically, we highlight that for $\alpha+\beta+1\neq0,$ Bernis'
inequality, see \cite[Lemma 2.1]{Ant}, yields, for $m\neq0$,
\begin{equation}\label{bernis conseq}
\frac19\int_{\mathbb T}
\frac{|\partial_x\rho^m|^4}{\rho^{2m}}\,dx
\le\int_{\mathbb T}|\partial_{xx}\rho^m|^2\,dx.
\end{equation}
In particular, the constant $\frac19$ is optimal, meaning that there
is no constant $C>\frac19$ such that \eqref{bernis conseq} holds for
all smooth positive periodic densities $\rho$.
By virtue of these considerations, we set
$m=\frac{\alpha+\beta+1}{2}$ and use \eqref{bernis conseq} to obtain
\begin{equation*}
\begin{split}
\frac{m^2}{c_\mu c_k}\mathcal F_k[\rho]
={}&\int_{\mathbb T}
\left|\partial_{xx}\rho^{\frac{\alpha+\beta+1}{2}}\right|^2\,dx+16\left[
\frac{(\alpha-\beta-1)(1-\alpha)}{(\alpha+\beta+1)^2}
-\frac{\beta}{3(\alpha+\beta+1)}
\right]
\int_{\mathbb T}
\left|\partial_x\rho^{\frac{\alpha+\beta+1}{4}}\right|^4\,dx\\
\ge{}&16\left[
\frac{(\alpha-\beta-1)(1-\alpha)}{(\alpha+\beta+1)^2}
-\frac{\beta}{3(\alpha+\beta+1)}+\frac19
\right]
\int_{\mathbb T}
\left|\partial_x\rho^{\frac{\alpha+\beta+1}{4}}\right|^4\,dx.
\end{split}
\end{equation*}
Since
\[
\frac{(\alpha-\beta-1)(1-\alpha)}{(\alpha+\beta+1)^2}
-\frac{\beta}{3(\alpha+\beta+1)}+\frac19
=\frac{2(2\alpha-\beta-1)(\beta-2\alpha+4)}
{9(\alpha+\beta+1)^2},
\]
this quantity is nonnegative for
$2\alpha-4\le\beta\le2\alpha-1.$ If $\alpha+\beta+1=0$, substitution of $\beta=-\alpha-1$ into the previous estimate yields 
\begin{equation*}
\mathcal F_k[\rho]
=c_\mu c_k\int_{\mathbb T}|\partial_{xx}\log\rho|^2\,dx
+\frac{c_\mu c_k\alpha(1-\alpha)}2
\int_{\mathbb T}|\partial_x\log\rho|^4\,dx.
\end{equation*}
The remaining terms are treated as in
Proposition \ref{Prop NS BD}, and the proof is concluded.
\end{proof}

\noindent
We conclude this section with the following density bounds which are
implied by Proposition \ref{Prop NSK Energy} and Proposition
\ref{Prop NSK BD}. We start with the following result which solely
relies on initial data in the energy space. At this level, the
friction term does not provide regularity gains.

\begin{proposition}
Let $(\rho,u)$ be a smooth solution of \eqref{NSK damping} with finite
initial energy. Then the following holds:
\begin{itemize}
\item[1)] For any $\beta\ge-2-\gamma$,
\begin{equation}\label{energy upper rho NSK}
\rho\in L^\infty((0,T)\times\mathbb T).
\end{equation}
\item[2)] For any $\beta\le-2$,
\begin{equation}\label{energy no vacuum NSK}
\frac1\rho\in L^\infty((0,T)\times\mathbb T).
\end{equation}
\end{itemize}
\end{proposition}

\begin{proof}
Similarly to the proof of Proposition \ref{Prop density NS},
we observe that for any $\beta\neq-2-\gamma$,
\begin{equation}\label{tfc upper rho energy NSK}
\begin{split}
\rho^{\frac{\beta+2+\gamma}{2}}(x,t)
&=1+\int_{x_0}^x
\partial_y\rho^{\frac{\beta+2+\gamma}{2}}\,dy \le1+\left|\frac{\beta+2+\gamma}{2}\right|
\int_{\mathbb T}
\rho^{\frac\gamma2}\rho^{\frac\beta2}|\partial_x\rho|\,dx
\le C.
\end{split}
\end{equation}
Therefore \eqref{energy upper rho NSK} follows by taking the
supremum in space and time in \eqref{tfc upper rho energy NSK},
using conservation of mass and Proposition \ref{Prop NSK Energy},
for $\beta>-2-\gamma$.
The proof for $\beta=-2-\gamma$ follows by using the same method
applied to $\log\rho$ in \eqref{tfc upper rho energy NSK}.
Analogously, we have
\begin{equation}\label{tfc lower rho energy NSK}
\rho^{\frac{\beta+2}{2}}(x,t)
\le1+\int_{\mathbb T}
\left|\partial_x\rho^{\frac{\beta+2}{2}}\right|\,dx.
\end{equation}
Thus for $\beta<-2$, \eqref{energy no vacuum NSK} holds.
Also in this case, the proof for $\beta=-2$ applies with minor
modifications by considering $\log\rho$ in the left-hand side
of \eqref{tfc lower rho energy NSK}.
\end{proof}

\noindent
The next proposition strongly uses the regularity inferred by the
BD entropy inequality \eqref{NSK BD ineq}. It represents the main
novel regularity estimate concerning the no-vacuum analysis of
compressible fluids with capillarity.

\begin{proposition}\label{Prop density NSK}
Let $(\rho,u)$ be a smooth solution of \eqref{NSK damping} with initial datum satisfying \eqref{reg initial datum}. Assume in addition that
\begin{equation}
2\alpha-4\le\beta\le2\alpha-1.
\end{equation}
Then the following upper and lower bounds on the density hold:
\begin{itemize}
\item[1)] For any $\alpha$ and $\beta$,
\begin{equation}\label{upper rho NSK}
\rho\in L^\infty((0,T)\times\mathbb T).
\end{equation}
\item[2)] For any $\alpha,\beta$ such that
$0\le\alpha\le\max\{\frac12,\frac{2-\theta}{3}\}$
or $\beta\le\max\{-2,-\alpha-\theta-1\}$,
\begin{equation}\label{no vacuum NSK}
\frac1\rho\in L^\infty((0,T)\times\mathbb T).
\end{equation}
\end{itemize}
\end{proposition}

\begin{remark}
For $\theta=0$, the bound \eqref{no vacuum NSK} holds whenever
$\alpha\in[0,\frac23]$ or $\beta\le-\alpha-1$, provided that
the strong coercivity condition \eqref{scc} holds.
These conditions reduce to \ref{epsilon NV NSK} and in particular
the capillarity exponent $\beta$ determines the absence of vacuum in the range
\begin{equation}
\frac23<\alpha\le1,
\qquad
2\alpha-4\le\beta\le-\alpha-1,
\end{equation}
and therefore values $\beta\in(-2,-\frac53)$ are allowed.
\end{remark}

\begin{proof}
Similarly to the proof of Proposition \ref{Prop density NS},
we recall the following inequality:
\begin{equation}\label{gen tfc rho}
\begin{split}
\rho^m(x,t)
&\le1+\int_{\mathbb T}|\partial_x\rho^m|\,dx=1+\left|\frac mn\right|
\int_{\mathbb T}\rho^{m-n}|\partial_x\rho^n|\,dx,
\qquad m,n\neq0,
\end{split}
\end{equation}
and
\begin{equation}\label{gen tfc log rho}
\begin{split}
|\log\rho(x,t)|
&\le\int_{\mathbb T}|\partial_x\log\rho|\,dx=\frac1{|n|}
\int_{\mathbb T}\rho^{-n}|\partial_x\rho^n|\,dx,
\qquad n\neq0.
\end{split}
\end{equation}
Then we observe that by virtue of Proposition \ref{Prop NSK Energy},
Proposition \ref{Prop NSK BD} and conservation of mass,
\eqref{upper rho NSK} follows by choosing $m=\alpha$ and
$n=\alpha-\frac12$ in \eqref{gen tfc rho}, when
$\alpha\notin\{0,\frac12\}$.
The cases $\alpha=0$ and $\alpha=\frac12$ are treated as in the proof of Proposition \ref{Prop density NS}.
The lower bound for $\alpha\le\frac12$ follows from the same argument, while the case $\beta\le-2$ is already covered by
\eqref{energy no vacuum NSK}. We may therefore assume
$\alpha>\frac12$ and $\beta>-2$.
On the other hand, by choosing again $n=\alpha-\frac12$
in \eqref{gen tfc log rho}, we also have
\begin{equation*}
\begin{split}
|\log\rho(x,t)|
&\le\frac1{|\alpha-\frac12|}
\int_{\mathbb T}
\rho^{\frac12-\alpha}
\left|\partial_x\rho^{\alpha-\frac12}\right|\,dx=\frac1{|\alpha-\frac12|}
\int_{\mathbb T}
\rho^{\frac{2-3\alpha-\theta}{2}}
\rho^{\frac{\alpha+\theta-1}{2}}
\left|\partial_x\rho^{\alpha-\frac12}\right|\,dx
\le C,
\end{split}
\end{equation*}
provided that $2-3\alpha-\theta\ge0$, which implies
$\alpha+\theta<1$.
Similarly, by using \eqref{gen tfc log rho} with
$n=\frac\beta2+1$, we have
\begin{equation*}
\begin{split}
|\log\rho(x,t)|
&\le\frac1{|\frac\beta2+1|}
\int_{\mathbb T}
\rho^{-\frac\beta2-1}
\left|\partial_x\rho^{\frac\beta2+1}\right|\,dx\\
&=\frac1{|\frac\beta2+1|}
\int_{\mathbb T}
\rho^{\frac{-\alpha-\beta-\theta-1}{2}}
\rho^{\frac{\alpha+\theta-1}{2}}
\left|\partial_x\rho^{\frac\beta2+1}\right|\,dx
\le C
\end{split}
\end{equation*}
whenever $\beta\le-\alpha-\theta-1$.
Since $\beta>-2$, this condition also implies $\alpha+\theta<1.$
Both inequalities include equality in the respective parameter
conditions. This concludes the proof of \eqref{no vacuum NSK}.
\end{proof}
\section{Pointwise estimates for the Navier-Stokes equations}\label{sec:active potential}
\noindent
Throughout this section, $(\rho,u)$ is a smooth positive solution of \eqref{NS damping} on $[0,T)\times\mathbb T$, with initial data satisfying \eqref{reg initial datum}. We write $m_0=\min_{\mathbb T}\rho_0>0$ and $M_0=\max_{\mathbb T}\rho_0$. We first use the effective pressure to obtain a lower density bound for a class of power-law friction terms. Then, we focus to the case of linear damping and we study the sign of the active potential. For sufficiently large damping, this also gives a quantitative bound on how far the density can move outside its initial range.

\subsection{Effective pressure}\label{subsec:effective-pressure}
We follow the effective-pressure argument of Burtea and Haspot~\cite[Section~3.1]{BurteaHaspot}. Recall that $V=u+Q$, where $Q=\frac{\mu(\rho)}{\rho^2}\partial_x\rho$. We seek an upper bound for $\frac{\partial_xV}{\rho}$ since, at a minimum of the density, this quantity controls how fast $\rho^{-1}$ can grow.

\begin{proposition}\label{Prop effective pressure NS}
Assume $\gamma\ge\alpha$, $f\in L^2(0,T;L^\infty(\mathbb T))$ and $\partial_xf\in L^1(0,T;L^\infty(\mathbb T))$. Suppose that $\partial_xV_0\le K_0\rho_0$ on $\mathbb T$ for some $K_0\ge0$, and that one of the following conditions holds.
\begin{itemize}
\item[(i)] $\theta=1$.
\item[(ii)] $\theta\ge\max\{1,\alpha\}$.
\item[(iii)] $\gamma>\alpha$ and $\theta\ge\frac{\gamma+1}{2}$.
\end{itemize}
Then
\begin{equation}\label{effective pressure density bound}
\sup_{0\le t<T}\left(
\left\|\frac{1}{\rho(t)}\right\|_{L^\infty}
+\left\|\left(\frac{\partial_xV(t)}{\rho(t)}\right)_+\right\|_{L^\infty}
\right)\le C(T).
\end{equation}
The constant depends only on $T$, $m_0$, $K_0$, $\varepsilon$, the constitutive parameters, the initial quantities in \eqref{reg initial datum} and the stated norms of $f$.
\end{proposition}

\begin{proof}
Write $D_t=\partial_t+u\partial_x$ and set
\begin{equation*}
A(r)=\frac{rp'(r)}{\mu(r)},\qquad b(r)=\varepsilon r^{\theta-1}.
\end{equation*}
The fluid velocity $u$ and the effective velocity $V=u+Q$ satisfy
\begin{equation}\label{effective velocity equations NS}
D_tu-\frac{1}{\rho}\partial_x\bigl(\mu(\rho)\partial_xu\bigr)
=D_tV=-A(\rho)Q-b(\rho)u+f.
\end{equation}
By the continuity equation, $D_t\bigl(\frac{\partial_xV}{\rho}\bigr)=\frac{1}{\rho}\partial_x(D_tV)$. To see how $\partial_xu$ enters this identity, we write $-A(\rho)Q-b(\rho)u=-A(\rho)V+(A(\rho)-b(\rho))u$. Differentiating in $x$ and dividing by $\rho$ then gives the coefficient $\frac{A(\rho)-b(\rho)}{\rho}$ in front of $\partial_xu$. This suggests introducing
\begin{equation*}
Y=\frac{\partial_xV}{\rho}+G(\rho),\qquad
G(r)=\int_1^r\frac{A(s)-b(s)}{s^2}\,ds.
\end{equation*}
Indeed, $D_tG(\rho)=-\rho G'(\rho)\partial_xu$ cancels this term, and the equation for the effective pressure becomes
\begin{equation}\label{effective pressure equation NS}
D_tY+A(\rho)Y
=A(\rho)G(\rho)
-\frac{\rho}{\mu(\rho)}\bigl(A'(\rho)Q^2+b'(\rho)uQ\bigr)
+\frac{\partial_xf}{\rho}.
\end{equation}
By \eqref{upper rho NS}, $\rho\le R_T$ for some $R_T\ge1$. Since $\gamma\ge\alpha$ and $\theta\ge1$ in all three cases, $A$ and $b$ are bounded on $(0,R_T]$, $A'\ge0$, and $|G(r)|\le C(T)(1+\frac{1}{r})$. The restrictions on $\theta$ arise in estimating the term containing $b'$. For $\theta=1$ this term vanishes. Under (ii), we use $Q=V-u$ and
\begin{equation*}
\frac{r|b'(r)|}{\mu(r)}
\le C r^{\theta-\alpha-1}
\le C(T)\left(1+\frac{1}{r}\right).
\end{equation*}
Under (iii), $A'>0$, so we can use the term in $Q^2$ to obtain
\begin{equation*}
-\frac{\rho}{\mu(\rho)}\bigl(A'(\rho)Q^2+b'(\rho)uQ\bigr)
\le\frac{\rho\,b'(\rho)^2}{4\mu(\rho)A'(\rho)}u^2
\le C\rho^{2\theta-\gamma-2}u^2
\le C(T)\left(1+\frac{1}{\rho}\right)u^2.
\end{equation*}
We now estimate $Y$ and $\rho^{-1}$ together. In terms of $V$, the continuity equation reads
\begin{equation*}
\partial_t\rho+\partial_x(\rho V)
=\partial_x\left(\frac{\mu(\rho)}{\rho}\partial_x\rho\right).
\end{equation*}
Set
\begin{equation*}
M(t)=\max\left\{\max_{x\in\mathbb T}Y(x,t),\,
                 \max_{x\in\mathbb T}\frac{1}{\rho(x,t)}\right\}
\end{equation*}
and choose $x_t\in\mathbb T$ where either $Y(\cdot,t)$ or $\rho(\cdot,t)^{-1}$ takes the value $M(t)$.
If $M(t)=\rho(x_t,t)^{-1}$, then $\rho(\cdot,t)$ attains its minimum at $x_t$. The diffusion term in the preceding equation is nonnegative there, and hence
\begin{equation*}
\partial_t(\rho^{-1})(x_t,t)
\le \frac{\partial_xV(x_t,t)}{\rho(x_t,t)}
=Y(x_t,t)-G(\rho(x_t,t)).
\end{equation*}
If $M(t)=Y(x_t,t)$, we use \eqref{effective pressure equation NS} at $x_t$. Since the mass is normalized to one, $M(t)\ge1$. Combining the two cases, we obtain for almost every $t$
\begin{equation}\label{eq:effective-pressure-gronwall}
M'\le C(T)\bigl(1+\|\partial_xf\|_{L^\infty}
+\|u\|_{L^\infty}^2+\|V\|_{L^\infty}^2\bigr)M.
\end{equation}
The bounds on $\rho_0$ and the assumption $\partial_xV_0\le K_0\rho_0$ give an upper bound for $M(0)$. When $\theta=1$, the terms involving the velocity norms drop out of \eqref{eq:effective-pressure-gronwall}, so Gronwall's inequality gives a bound for $M$ on $[0,T]$. Using $\frac{\partial_xV}{\rho}=Y-G(\rho)$, we obtain \eqref{effective pressure density bound} and conclude the proof for $\theta=1$. 
\\
\\
For $\theta>1$, we first bound $V_0$. Since $\partial_xV_0\le K_0\rho_0$ and the initial mass is equal to one, we have
\begin{equation*}
\max_{\mathbb T}V_0-\min_{\mathbb T}V_0
\le \int_{\mathbb T}(\partial_xV_0)_+\,dx\le K_0,
\end{equation*}
from which we also deduce $\|V_0\|_{L^\infty}\le K_0+\|\sqrt{\rho_0}V_0\|_{L^2}.$ Next, let $\tau$ be the first time at which $\min_{\mathbb T}\rho=\frac{m_0}{2}$, with $\tau=T$ if this value is never reached. On $(0,\tau)$ we have $\rho\ge\frac{m_0}{2}$, so the energy estimate \eqref{NS energy ineq} and the embedding $H^1(\mathbb T)\subset L^\infty(\mathbb T)$ give
\begin{equation*}
\int_0^\tau\|u(t)\|_{L^\infty}^2\,dt\le C(T).
\end{equation*}
Since $D_tV+A(\rho)V=(A(\rho)-b(\rho))u+f$, with $A\ge0$ and $A,b$ bounded, we also have, for $t<\tau$,
\begin{equation*}
\|V(t)\|_{L^\infty}\le \|V_0\|_{L^\infty}
+\int_0^t\bigl(C(T)\|u(r)\|_{L^\infty}+\|f(r)\|_{L^\infty}\bigr)\,dr\le C(T).
\end{equation*}
Gronwall's inequality in \eqref{eq:effective-pressure-gronwall} now bounds $M$ on $[0,\tau)$. If $\tau=T$, this gives \eqref{effective pressure density bound}.
\\
Suppose now that $\tau<T$, then the estimate at a minimum of $\rho$ gives
\begin{equation*}
\frac{\mathrm d}{\mathrm dt}\left(\max_{x\in\mathbb T}\frac{1}{\rho(x,t)}\right)
\le C(T)
\qquad\text{for almost every }t\in(0,\tau).
\end{equation*}
At time $\tau$ the minimum density is $\frac{m_0}{2}$, so integrating from $0$ to $\tau$ gives
\begin{equation*}
\frac{1}{m_0}
=\max_{x\in\mathbb T}\frac{1}{\rho(x,\tau)}
-\max_{x\in\mathbb T}\frac{1}{\rho_0(x)}
\le C(T)\tau.
\end{equation*}We can thus choose $s\in(0,\tau)$ such that
\begin{equation*}
\|u(s)\|_{L^\infty}^2\le\frac{1}{\tau}\int_0^\tau\|u(t)\|_{L^\infty}^2\,dt\le C(T).
\end{equation*}
At this time we already have $\|V(s)\|_{L^\infty}\le C(T)$. Starting from $s$, the maximum principle for $u$ and the transport equation for $V$ give
\begin{equation*}
\frac{d}{dt}\bigl(\|u\|_{L^\infty}+\|V\|_{L^\infty}\bigr)
\le C(T)\bigl(\|u\|_{L^\infty}+\|V\|_{L^\infty}\bigr)+2\|f\|_{L^\infty},
\end{equation*}
since both equations in \eqref{effective velocity equations NS} have right-hand side $(A-b)u-AV+f$. Gronwall's inequality bounds both velocities on $[s,T)$, and \eqref{eq:effective-pressure-gronwall} then gives \eqref{effective pressure density bound} on the whole interval.
\end{proof}
\subsection{Active potential with linear damping}\label{subsec:active-potential-linear-damping}

We now take $\theta=1$ and a spatially homogeneous force $f=f(t)\in L^2(0,T)$. Following \cite{Const}, we introduce the active potential $w=\mu(\rho)\partial_xu-p(\rho)$. The velocity equation then reads
\begin{equation}\label{NS linear damping}
D_tu=\frac{1}{\rho}\partial_xw-\varepsilon u+f(t).
\end{equation}
Differentiating \eqref{NS linear damping} in space and using the continuity equation gives
\begin{equation}\label{eq:w-damped}
\begin{split}
\partial_tw={}&c_\mu\rho^{\alpha-1}\partial_{xx}w
-\bigl(u+c_\mu\rho^{\alpha-2}\partial_x\rho\bigr)\partial_xw\\
&-(w+p(\rho))\left[\varepsilon+
\frac{(\alpha+1)w+(\alpha+1-\gamma)p(\rho)}{\mu(\rho)}\right].
\end{split}
\end{equation}
Write $\rho_m,\rho_M$ and $w_m,w_M$ for the spatial minima and maxima, and set $W_0=(w_M(0))_+$. These extrema are locally absolutely continuous. At almost every time, the reaction term in \eqref{eq:w-damped}, evaluated at a maximum point of $w$, gives an upper bound for $w_M'$. At a minimum point it gives a lower bound for $w_m'$. At the density extrema the continuity equation gives
\begin{equation}\label{eq:ap-density-minimum}
\rho_m'\ge-\frac{\rho_m}{\mu(\rho_m)}\bigl(w_M+p(\rho_m)\bigr),
\qquad
\rho_M'\le-\frac{\rho_M}{\mu(\rho_M)}\bigl(w_m+p(\rho_M)\bigr).
\end{equation}
\begin{proposition}\label{prop:ap-sign-entry}
If $\gamma\le\alpha+1$, there exists a time $T_\varepsilon\ge0$, depending only on $\varepsilon$, $W_0$, the initial bounds in \eqref{reg initial datum} and the constitutive parameters, such that whenever $T>T_\varepsilon$,
\begin{equation}\label{eq:ap-first-entry}
w(x,t)\le0 \qquad \text{for all } (x,t)\in\mathbb T\times[T_\varepsilon,T).
\end{equation}
Assume that one of the following conditions holds
\begin{itemize}
\item[(i)] $\gamma\ge\alpha$ and $\partial_xV_0\le K_0\rho_0$ on $\mathbb T$ for some $K_0\ge0$.
\item[(ii)] $\alpha\le\gamma\le\alpha+1$ and $w_0\le0$ on $\mathbb T$.
\end{itemize}
Then there exists a constant $c>0$, independent of $T$, such that
\begin{equation*}
\rho(x,t)\ge c \qquad\text{for all }(x,t)\in\mathbb T\times[0,T).
\end{equation*}
\end{proposition}

\begin{proof}
We may assume without loss of generality that $f=0$. Indeed, let $X$ solve $X''+\varepsilon X'=f$ with $X(0)=X'(0)=0$. Then $\rho(x+X(t),t)$ and $u(x+X(t),t)-X'(t)$ solve the unforced problem with the same initial data, while this translation preserves the bounds on $\rho$ and $w$.
\\
Then, the main observation is that the total time during which the minimum density lies below half its mean is bounded. The unforced BD balance in the proof of Proposition~\ref{Prop NS BD} gives
\begin{equation*}
\int_0^T\int_{\mathbb T}\frac{\mu(\rho)p'(\rho)}{\rho^2}|\partial_x\rho|^2\,dx\,dt
\le E_0,\qquad
E_0=H^{BD}(\rho_0,u_0)+\int_{\mathbb T}I(\rho_0;\varepsilon)\,dx.
\end{equation*}
Here $I\ge0$, and $\int_{\mathbb T}h(\rho)\,dx\ge h(1)=0$, since the mean density is one. In particular, $\rho_M(t)\ge1$. Thus, whenever $\rho_m(t)<\frac12$, the density takes every value between $\frac12$ and $1$. By the fundamental theorem of calculus and the Cauchy-Schwarz inequality,
\begin{equation*}
0<\delta:=\int_{1/2}^1\frac{\sqrt{\mu(r)p'(r)}}{r}\,dr
\le \int_{\mathbb T}\frac{\sqrt{\mu(\rho)p'(\rho)}}{\rho}|\partial_x\rho|\,dx
\le \left(\int_{\mathbb T}\frac{\mu(\rho)p'(\rho)}{\rho^2}|\partial_x\rho|^2\,dx\right)^{1/2}.
\end{equation*}
The pressure dissipation is therefore at least $\delta^2$ at each such time. Integrating in time and using the preceding BD estimate, we obtain
\begin{equation}\label{eq:ap-low-density-times}
E=\left\{t\in[0,T):\rho_m(t)<\tfrac12\right\},
\qquad |E|\le B:=\frac{E_0}{\delta^2}.
\end{equation}
We emphasize that the choice of $\frac{1}{2}$ is only for convenience of presentation, since the same argument applies to any fixed threshold $a\in(0,1)$.  Suppose first that $\gamma\le\alpha+1,$ then at $w=0$ the reaction in \eqref{eq:w-damped} is strictly negative, so the maximum principle preserves $w\le0$. While $w_M>0$, the same equation gives
\begin{equation}\label{eq:ap-positive-decay}
w_M'\le-\varepsilon\bigl(w_M+p(\rho_m)\bigr).
\end{equation}
It follows that $\log\left(1+\frac{w_M}{p(\frac12)}\right)$ is nonincreasing while $w_M>0$, and its derivative is at most $-\varepsilon$ whenever $t\notin E$. Thus, for $W_0>0$, we may take
\begin{equation}\label{eq:ap-entropy-entry-time}
T_\varepsilon=B+\frac{1}{\varepsilon}\log\left(1+\frac{W_0}{p(\frac12)}\right).
\end{equation}
Suppose $T>T_\varepsilon$. Since $|E|\le B$, we have $|[0,T_\varepsilon]\setminus E|\ge T_\varepsilon-B$. If $w_M$ remained positive on $[0,T_\varepsilon]$, integration would give
\begin{equation*}
\log\left(1+\frac{w_M(T_\varepsilon)}{p(\frac12)}\right)
\le \log\left(1+\frac{W_0}{p(\frac12)}\right)-\varepsilon(T_\varepsilon-B)=0,
\end{equation*}
contradicting $w_M(T_\varepsilon)>0$. Once $w_M$ reaches zero, the maximum principle gives \eqref{eq:ap-first-entry}.
\\
We now prove the lower density bound under assumption (ii). In this case $w\le0$ from the initial time, and \eqref{eq:ap-density-minimum} gives
\begin{equation}\label{eq:ap-first-density}
\rho_m'\ge-\frac{c_p}{c_\mu}\rho_m^{1+\gamma-\alpha}.
\end{equation}
On the set $E,$ we have $\rho_m<\frac12$, so $\rho_m^{\gamma-\alpha}\le1$ since $\gamma\ge\alpha$. Dividing \eqref{eq:ap-first-density} by $\rho_m$ therefore gives $(\log\rho_m)'\ge-\frac{c_p}{c_\mu}$ on $E$. Outside $E$, the derivative of $\log\min\{\rho_m,\frac12\}$ vanishes almost everywhere. Integrating and using $|E|\le B$, we obtain
\begin{equation*}
\rho_m(t)\ge\min\left\{m_0,\tfrac12\right\}
\exp\left(-\frac{c_p}{c_\mu}B\right).
\end{equation*}
\noindent
If $\partial_xV_0\le K_0\rho_0$ and $\alpha\le\gamma\le\alpha+1$, Proposition~\ref{Prop effective pressure NS} gives a lower bound for the density up to $T_\varepsilon$. From $T_\varepsilon$ onwards, $w\le0$ and the preceding argument applies. To obtain the density bound under assumption (i) for the full range $\gamma\ge\alpha$, we return to the effective pressure. In the unforced system, the energy and the modified BD entropy are nonincreasing. The proof of \eqref{upper rho NS} therefore gives an upper bound $\rho\le R$ independent of $T$.
\\
\\
With $Y$ and $G$ as in the preceding subsection, choose $K>2$ large enough that $G(r)\le K$ for every $r\in[\frac12,R]$, and set
\begin{equation*}
M_K(t)=\max\left\{K,\max_{\mathbb T}Y(t),\frac{1}{\rho_m(t)}\right\}.
\end{equation*}
When $t\notin E$, we have $\rho_m(t)\ge\frac12$, so $\frac{1}{\rho_m(t)}\le2$ and $G(\rho)\le K$. At a spatial maximum of $Y$ with $Y\ge K$, equation \eqref{effective pressure equation NS} gives
\begin{equation*}
D_tY\le A(\rho)\bigl(G(\rho)-Y\bigr)\le0,
\end{equation*}
since $b'=0$, $\partial_xf=0$ and $A'\ge0$. Thus $M_K'\le0$ almost everywhere outside $E$. On $E$, the estimate leading to \eqref{eq:effective-pressure-gronwall} gives $M_K'\le C M_K$. Here $\theta=1$, so the velocity norms do not enter the estimate, and $C$ depends only on $R$ and the coefficients. Since $|E|\le B$, Gronwall's inequality gives
\begin{equation*}
\rho_m(t)\ge \frac{1}{M_K(t)}\ge \frac{e^{-CB}}{M_K(0)}.
\end{equation*}
The bounds on $\rho_0$ and the assumption $\partial_xV_0\le K_0\rho_0$ control $M_K(0)$. This completes the proof.
\end{proof}
\noindent
The uniform estimates from below for the density can also be obtained without damping under either hypothesis (i) or (ii) of Proposition~\ref{prop:ap-sign-entry}. When $\gamma<\alpha+1$, the argument for the sign of $w$ also applies in the undamped case. Indeed, while $w_M>0$, the pressure term alone makes $w_M$ decrease at a rate bounded away from zero whenever $\rho_m\ge\frac12$. Together with \eqref{eq:ap-low-density-times}, this gives a different estimate for the time after which $w$ remains nonpositive, provided the solution exists that long. The role of damping is particularly clear at $\gamma=\alpha+1$, where the pressure contribution vanishes at $w=0$ and the reaction reduces to $-\varepsilon p(\rho)<0$.
\\
\\
We finally consider arbitrary $\alpha\ge0$ and $\gamma\ge1$. When $\gamma<\alpha$, even if $w\le0$, integrating \eqref{eq:ap-first-density} does not give a positive lower bound for the density for all times. For sufficiently large $\varepsilon$, we can instead estimate how far the minimum density can fall below $m_0$ and how far the maximum density can rise above $M_0$.

\begin{proposition}\label{active potential extended thm}
Let $\alpha\ge0$ and $\gamma\ge1$. There exist $C,\varepsilon_0>0$, depending only on $m_0$, $M_0$, $\|w_0\|_{L^\infty}$ and the constitutive parameters, such that for every $\varepsilon\ge\varepsilon_0$,
\begin{equation}\label{eq:ap-density-range}
m_0-\frac{C}{\varepsilon}\le\rho(x,t)\le M_0+\frac{C}{\varepsilon}
\qquad\text{on }\mathbb T\times[0,T).
\end{equation}
Furthermore, if $0<L<m_0\le M_0<R$, we may choose $\varepsilon_0$ sufficiently large, depending also on $L$ and $R$, so that for every $\varepsilon\ge\varepsilon_0$ the density stays in $[L,R]$ and the following estimates hold for all $0\le t<T$:
\begin{equation}\label{eq:ap-second-potential-bounds}
\begin{gathered}
\min\{w_m(0),-p(R)\}\le w(x,t)\le\max\{w_M(0),-p(L)\},\\
\|\partial_xu(t)\|_{L^\infty}\le\frac{\|w_0\|_{L^\infty}+p(R)}{\mu(L)}.
\end{gathered}
\end{equation}
Moreover, there exists $C_{L,R}>0$, independent of $T$ and $\varepsilon$, such that for $\frac{C_{L,R}}{\varepsilon}\le t<T$,
\begin{equation}\label{eq:ap-second-negative}
-p(R)\le w(x,t)\le-p(L),
\qquad
\|\partial_xu(t)\|_{L^\infty}\le\frac{p(R)-p(L)}{\mu(L)}.
\end{equation}
\end{proposition}

\begin{proof}
We briefly describe the idea behind the proof of Proposition~\ref{active potential extended thm}. Given $0<L<m_0\le M_0<R$, we show that, for sufficiently large damping, the rectangle $[L,R]\times[-p(R),-p(L)]$ is invariant for $(\rho,w)$. In other words, if the bounds $L\le\rho\le R$ and $-p(R)\le w\le-p(L)$ hold throughout $\mathbb T$ at some time, then both continue to hold for as long as the solution exists. \\
The threshold $\varepsilon_0$ can be chosen explicitly and serves two purposes. First, since $w_0$ may lie outside $[-p(R),-p(L)]$, we choose the damping large enough that the reaction term in \eqref{eq:w-damped} drives $w$ toward this interval as long as $L\le\rho\le R$. Second, we need to control the density before $w$ enters the interval. The exponential bounds following \eqref{eq:ap-barrier-maxima} give the one-sided estimates \eqref{eq:ap-density-range-one-sided}: the density can fall below $m_0$ or rise above $M_0$ by at most $\frac{C}{\varepsilon}$. By increasing $\varepsilon_0$, we make this distance smaller than both $m_0-L$ and $R-M_0$. The density then stays strictly between $L$ and $R$, and the same exponential bounds give \eqref{eq:ap-second-negative} by time $\frac{C_{L,R}}{\varepsilon}$, provided the solution exists that long.
\\
\\
Set $J=[\frac{m_0}{2},2M_0]$ and consider an interval $[0,t]$ on which $\rho$ takes values in $J$. Let $w_*=\min\{w_m(0),-p(2M_0)\}$. Since the expression below is increasing in $w$ and $J$ is compact, we may choose $\varepsilon_0$ so that
\begin{equation}\label{eq:ap-barrier-maxima}
\varepsilon+\frac{(\alpha+1)w+(\alpha+1-\gamma)p(s)}{\mu(s)}
\ge\frac{\varepsilon}{2}
\qquad\text{for }s\in J,\quad w\ge w_*.
\end{equation}
At $w=w_*$ we have $w_*+p(\rho)\le0$, so the reaction in \eqref{eq:w-damped} is nonnegative. The maximum principle gives $w\ge w_*$ throughout $[0,t]$.
\\
\\
Next, suppose that $a\le\rho\le b$ on $[0,t]$, for some $a,b>0$. At a maximum of $w$, whenever $w_M+p(a)>0$, the reaction is at most $-\frac{\varepsilon}{2}(w_M+p(a))$. At a minimum, whenever $w_m+p(b)<0$, it is at least $-\frac{\varepsilon}{2}(w_m+p(b))$. Integrating the corresponding inequalities gives
\begin{equation*}
\begin{aligned}
\bigl(w_M(t)+p(a)\bigr)_+ &\le \bigl(w_M(0)+p(a)\bigr)_+e^{-\varepsilon t/2},\\
\bigl(-w_m(t)-p(b)\bigr)_+ &\le \bigl(-w_m(0)-p(b)\bigr)_+e^{-\varepsilon t/2}.
\end{aligned}
\end{equation*}
\noindent
For each fixed $t$, we can apply this comparison with $a=\ell(t)$ and $b=r(t)$, where
\begin{equation*}
\ell(t)=\min_{0\le s\le t}\rho_m(s),\qquad
r(t)=\max_{0\le s\le t}\rho_M(s).
\end{equation*}
The function $\ell$ can decrease only when $\ell=\rho_m$, and $r$ can increase only when $r=\rho_M$. Thus \eqref{eq:ap-density-minimum} gives, almost everywhere,
\begin{equation*}
-\ell'\le a_0\bigl(w_M+p(\ell)\bigr)_+,\qquad
r'\le a_0\bigl(-w_m-p(r)\bigr)_+,\qquad
a_0=\max_{s\in J}\frac{s}{\mu(s)}.
\end{equation*}
Since $\ell(t)\le m_0$ and $r(t)\ge M_0$, the exponential estimates above give, after integration,
\begin{equation}\label{eq:ap-density-range-one-sided}
\begin{aligned}
\ell(t)&\ge m_0-\frac{2a_0}{\varepsilon}\bigl(w_M(0)+p(m_0)\bigr)_+,\\
r(t)&\le M_0+\frac{2a_0}{\varepsilon}\bigl(-w_m(0)-p(M_0)\bigr)_+.
\end{aligned}
\end{equation}
For sufficiently large $\varepsilon$, the lower bound is strictly larger than $\frac{m_0}{2}$ and the upper bound is strictly smaller than $2M_0$. By continuity, $\rho$ stays in $J$ throughout $[0,T)$, and \eqref{eq:ap-density-range} follows.
\\
\\
Finally, fix $L$ and $R$ as in the statement. Since $\mu$ is positive on $[L,R]$, we may increase $\varepsilon_0$ so that the bracket in \eqref{eq:w-damped} is nonnegative throughout the rectangle $[L,R]\times[-p(R),-p(L)]$. The reaction is then nonnegative at $w=-p(R)$ and nonpositive at $w=-p(L)$. In turn, whenever $-p(R)\le w\le-p(L)$, the continuity equation gives $D_t\rho\ge0$ at $\rho=L$ and $D_t\rho\le0$ at $\rho=R$. This proves the invariance of the rectangle. To show that the solution enters it, we use \eqref{eq:ap-density-range} and increase $\varepsilon_0$ further so that
\begin{equation*}
a:=\frac{L+m_0}{2}\le\rho(x,t)\le\frac{R+M_0}{2}=:b
\qquad\text{on }\mathbb T\times[0,T).
\end{equation*}
With these choices of $a$ and $b$, the same exponential estimates give the bounds on $w$ in \eqref{eq:ap-second-potential-bounds} for all $t$. Once $\varepsilon t$ is sufficiently large, the two right-hand sides are at most $p(a)-p(L)$ and $p(R)-p(b)$, respectively, since both differences are strictly positive. This gives the bounds on $w$ in \eqref{eq:ap-second-negative}, while the two estimates for $\partial_xu$ now follow from
\begin{equation*}
\partial_xu=\frac{w+p(\rho)}{\mu(\rho)},\qquad
p(L)\le p(\rho)\le p(R),\qquad \mu(\rho)\ge\mu(L),
\end{equation*}
this concludes the proof.
\end{proof}
\begin{remark}\label{remark high-friction}
Proposition \ref{active potential extended thm} can also be interpreted in light of the high-friction relaxation limit.  In the unforced case, the formal limit of the Navier-Stokes equations with damping $-\varepsilon\rho u$, as $\varepsilon\to\infty$ on the time scale $\tau=\frac{t}{\varepsilon}$, gives
\begin{equation*}
\partial_\tau\bar\rho=\partial_{xx}p(\bar\rho),
\qquad \bar\rho(x,0)=\rho_0(x).
\end{equation*}
This is the porous medium equation for $\gamma>1$ and the heat equation for $\gamma=1$.  By the maximum principle, comparison with the constant solutions $m_0$ and $M_0$ gives $m_0\le\bar\rho(x,\tau)\le M_0$ for all $\tau\ge0$. For the Navier-Stokes equations at finite friction, the initial velocity may drive the density outside $[m_0,M_0]$. Proposition~\ref{active potential extended thm} shows that, once $\varepsilon$ exceeds a threshold depending on the initial data and the constitutive parameters, the density can fall below $m_0$ or rise above $M_0$ by at most $\frac{C}{\varepsilon}$, uniformly in time. This estimate requires no smallness assumption on the initial data.
\end{remark}
\section{Global well-posedness}\label{Sec global wp}
\noindent
We devote this section to the proof of two global well-posedness results for the compressible Navier-Stokes and Navier-Stokes-Korteweg equations, stated in Theorem \ref{thm:NS-global} and Theorem \ref{thm:NSK-global}. The proofs follow a standard continuation argument. We first recall some local well-posedness results and then show how the no-vacuum estimates established in the previous sections allow us to extend the local solutions globally in time. For the Korteweg system, the BD entropy plays a key role in deriving the higher-order estimates needed for continuation. We only sketch the proofs, focusing on the propagation of the required regularity. 
\\
\\
We suggest comparing Theorem \ref{thm:NS-global} with the undamped Navier-Stokes results on the torus \cite[Theorem 1.6]{Const} and on $\mathbb R$ \cite[Theorem 2.1]{Mellet}, \cite[Theorem 2.2]{BurteaHaspot}. For Korteweg fluids, a direct comparison is given by the unforced deterministic case of \cite[Theorem 2.4]{Pesc}. See also \cite[Theorem~4.1]{Ant} for global existence and uniqueness of strong solutions to an approximating system with regularized viscosity and capillarity coefficients.

\subsection{Navier-Stokes equations}
This subsection contains the global well-posedness result for the Navier-Stokes equations in the improved no-vacuum regimes obtained in Sections \ref{Sec2} and \ref{sec:active potential}. We start by recalling the local well-posedness class.
\\
\\
Let $k\ge3$ be an integer and assume that $(\rho_0,u_0)\in H^k(\mathbb T)\times H^k(\mathbb T)$ with $\min_{\mathbb T}\rho_0>0$. Let $f\in L^2_{\mathrm{loc}}([0,\infty);H^{k-1}(\mathbb T))$. Proposition B.1 in \cite{Const} provides local existence and uniqueness for the system without friction. The friction term can be included in this construction, since it is a lower order term and the coefficient $\rho\mapsto\varepsilon\rho^{\theta-1}$ is smooth away from zero. Thus there exists a unique maximal local strong solution $(\rho,u)$ to system \eqref{NS damping}, with maximal existence time $T^*\in(0,\infty]$, such that
\begin{equation}\label{eq:NS-local-regularity}
\begin{split}
\rho,u&\in C([0,T^*);H^k(\mathbb T)),\quad
u\in L^2_{\mathrm{loc}}([0,T^*);H^{k+1}(\mathbb T)),
\end{split}
\end{equation}
and $\rho(x,t)>0$ on $\mathbb T\times[0,T^*)$. The solution can be continued as long as $\|\rho^{-1}\|_{L^\infty}$, $\|\rho\|_{H^k}$ and $\|u\|_{H^k}$ remain bounded. The following theorem collects the different conditions under which the preceding density estimates imply $T^*=\infty$.

\begin{theorem}\label{thm:NS-global}
Let $\alpha,\theta\ge0$, $\gamma\ge1$, $\varepsilon>0$, and let $k\ge3$ be an integer. Assume that
\begin{equation}\label{eq:NS-global-data}
(\rho_0,u_0)\in H^k(\mathbb T)\times H^k(\mathbb T),
\qquad
m_0=\min_{\mathbb T}\rho_0>0,
\end{equation}
and
\begin{equation*}
f\in L^2_{\mathrm{loc}}\bigl([0,\infty);H^{k-1}(\mathbb T)\bigr).
\end{equation*}
Suppose that one of the following conditions holds.
\begin{itemize}
\item[(i)]
$0\le\alpha\le\max\{\frac12,\frac{2-\theta}{3}\}$.
\item[(ii)]
$\gamma\ge\alpha$ and at least one of the following holds:
\par
{\centering
$\displaystyle \theta=1,\qquad \theta\ge\max\{1,\alpha\},\qquad
\text{or}\qquad \gamma>\alpha \ \text{and}\ \theta\ge\frac{\gamma+1}{2}.$
\par}
\item[(iii)]
$\theta=1$, $f(x,t)=f(t)$ and $\varepsilon\ge\varepsilon_0$.
\end{itemize}
Here the threshold $\varepsilon_0>0$ in (iii) is given in Proposition \ref{active potential extended thm} and depends only on $m_0$, $M_0$, $\|w_0\|_{L^\infty}$ and the constitutive parameters.
Then system \eqref{NS damping} admits a unique global strong solution satisfying, for every $T>0$,
\begin{equation}\label{eq:NS-global-regularity}
\begin{split}
\rho,u&\in C([0,T];H^k(\mathbb T)),\\
u&\in L^2(0,T;H^{k+1}(\mathbb T)),
\end{split}
\end{equation}
and
\begin{equation}\label{eq:NS-global-density}
0<c_T\le\rho(x,t)\le C_T
\qquad\text{on }\mathbb T\times[0,T].
\end{equation}
For $\theta=1$ and $f(x,t)=f(t)$, the bounds in \eqref{eq:NS-global-density} are uniform in time whenever $\gamma\ge\alpha$ or $\varepsilon\ge\varepsilon_0$.
\end{theorem}
\begin{remark}
Some comments are in order:
\begin{itemize}
\item[(1)] For $\varepsilon=0$, Burtea and Haspot \cite{BurteaHaspot} prove global existence and uniqueness of strong solutions on the real line $\mathbb R$ for $\alpha\ge\frac12$ and $\gamma\ge\max\{1,\alpha\}$. Our adaptation of the effective-pressure argument applies under the assumptions of Theorem~\ref{thm:NS-global}(ii), which require $\theta\ge1$. For $0\le\theta<\frac12$, the modified BD entropy gives a different improvement: for every $\varepsilon>0$, it raises the no-vacuum threshold of the undamped BD estimate from $\alpha\le\frac12$ to $\alpha\le\frac{2-\theta}{3}$, without imposing a one-sided bound on $\partial_xV_0$.
\\
\item[(2)] For linear damping $\theta=1$ and $f(x,t)=f(t)$, by virtue of the active-potential argument, we prove global existence and uniqueness of strong solutions for every $\alpha\ge0$ and $\gamma\ge1$, provided that $\varepsilon$ is sufficiently large. In particular, this covers the range $\gamma<\alpha$, which remains open in general for the undamped periodic problem with positive initial density. 
\end{itemize}
\end{remark}
\begin{proof}
Let $(\rho,u)$ be the maximal local strong solution in the class \eqref{eq:NS-local-regularity}, and fix $T<T^*$. We sketch the continuation argument, following the standard approach in \cite{Mellet,Const}. The estimates of Section \ref{sec:active potential} apply to the local strong solution, with the differential inequalities for the spatial extrema understood almost everywhere in time.
\\
In case (i), Proposition \ref{Prop density NS} gives \eqref{eq:NS-global-density}. In case (ii), the assumptions on the initial data give $V_0\in H^{k-1}(\mathbb T)$. Since $k\ge3$ and $\rho_0\ge m_0>0$, it follows that $\frac{\partial_xV_0}{\rho_0}\in L^\infty(\mathbb T)$. The assumption on $f$ also gives $\partial_xf\in L^1(0,T;L^\infty(\mathbb T))$. Thus the lower bound follows from Proposition \ref{Prop effective pressure NS}, while the upper bound follows from Proposition \ref{Prop density NS}. In case (iii), we apply Proposition \ref{active potential extended thm} with $L=\frac{m_0}{2}$ and $R=2M_0$. This gives \eqref{eq:NS-global-density} with bounds independent of $T$, as well as a uniform bound for $\|\partial_xu(t)\|_{L^\infty}.$
\\
\\
The energy and BD entropy estimates, together with the density bounds, give
\begin{equation}\label{eq:NS-global-base}
\sup_{t\in[0,T]}
\left(\|u(t)\|_{L^2}^2+\|\partial_x\rho(t)\|_{L^2}^2\right)
+\int_0^T\|\partial_xu\|_{L^2}^2\,dt\le C_T.
\end{equation}
We write the velocity equation as
\begin{equation}\label{eq:NS-global-velocity}
\begin{split}
\partial_tu+u\partial_xu
-\frac{\mu(\rho)}{\rho}\partial_{xx}u
&=\frac{\mu'(\rho)}{\rho}\partial_x\rho\,\partial_xu
-\frac{p'(\rho)}{\rho}\partial_x\rho-\varepsilon\rho^{\theta-1}u+f.
\end{split}
\end{equation}
Its diffusion coefficient is bounded away from zero. Testing \eqref{eq:NS-global-velocity} with $-\partial_{xx}u$ and using the H\"older and Sobolev inequalities, we have
\begin{equation*}
\frac{d}{dt}\|\partial_xu\|_{L^2}^2
+c_T\|\partial_{xx}u\|_{L^2}^2
\le C_T\left(1+\|\partial_xu\|_{L^2}^2+\|f\|_{L^2}^2\right)
\left(1+\|\partial_xu\|_{L^2}^2\right).
\end{equation*}
Since $\|\partial_xu\|_{L^2}^2\in L^1(0,T)$ by \eqref{eq:NS-global-base}, Gronwall's lemma gives
\begin{equation}\label{eq:NS-global-H1}
u\in L^\infty(0,T;H^1)\cap L^2(0,T;H^2),
\qquad
\partial_xu\in L^2(0,T;L^\infty).
\end{equation}
We next estimate higher derivatives. In particular, we define
\begin{equation*}
E_j=\|\rho\|_{H^j}^2+\|u\|_{H^j}^2.
\end{equation*}
At order $j\ge2$, we differentiate the continuity equation $j$ times and test with $\partial_x^j\rho$. For the velocity equation, we apply only $\partial_x^{j-1}$ and test with $-\partial_x^{j+1}u$. This gives the $H^j$ estimate for $u$ without introducing a derivative of order $j+1$ on the density. The term containing $\rho\partial_x^{j+1}u$ in the differentiated continuity equation is absorbed by the viscous dissipation through Young's inequality.
\\
For $j=2$, the bounds already obtained and one-dimensional interpolation give
\begin{equation*}
\|\partial_x\rho\|_{L^\infty}
\le C_T\left(1+\|\partial_{xx}\rho\|_{L^2}^{1/2}\right),
\qquad
\|\partial_{xx}u\|_{L^2}
\le C_T\|\partial_{xxx}u\|_{L^2}^{1/2}.
\end{equation*}
For example, for every $\delta>0$,
\begin{equation*}
\|\partial_x\rho\|_{L^\infty}
\|\partial_{xx}u\|_{L^2}
\|\partial_{xxx}u\|_{L^2}
\le\delta\|\partial_{xxx}u\|_{L^2}^2
+C_{\delta,T}\left(1+\|\partial_{xx}\rho\|_{L^2}^2\right).
\end{equation*}
Estimating the remaining terms in the same way yields
\begin{equation*}
\frac{d}{dt}E_2+c_T\|u\|_{H^3}^2
\le C_T\left(1+\|\partial_xu\|_{L^\infty}^2\right)(1+E_2)
+C_T\|f\|_{H^1}^2.
\end{equation*}
Since $1+\|\partial_xu\|_{L^\infty}^2$ is integrable in time by \eqref{eq:NS-global-H1}, Gronwall's lemma gives
\begin{equation}\label{eq:NS-global-H2}
\rho,u\in L^\infty(0,T;H^2),
\qquad
u\in L^2(0,T;H^3).
\end{equation}
For each $3\le j\le k$, the same differentiated estimates and the bounds in \eqref{eq:NS-global-H2} give
\begin{equation*}
\frac{d}{dt}E_j+c_T\|u\|_{H^{j+1}}^2
\le C_T\left(1+\|u\|_{H^3}^2\right)(1+E_j)
+C_T\|f\|_{H^{j-1}}^2.
\end{equation*}
The constants here also depend on the bounds in \eqref{eq:NS-global-H2}. A further application of Gronwall's lemma therefore gives
\begin{equation*}
\sup_{t\in[0,T]}E_k(t)
+\int_0^T\|u(t)\|_{H^{k+1}}^2\,dt\le C_T.
\end{equation*}
The continuity in time in \eqref{eq:NS-global-regularity} follows from the local theory and the continuity equation.
\\
If $T^*<\infty$, the density and Sobolev bounds remain uniform as $T\uparrow T^*$. The local continuation criterion then allows us to extend the solution beyond $T^*$, which is a contradiction. Thus $T^*=\infty$. Finally, uniqueness follows from the local strong solution theory on successive time intervals.
\end{proof}

\subsection{Navier-Stokes-Korteweg equations}

We now recall the local theory for system \eqref{NSK damping}.
Let $k\ge3$ be an integer, assume that $(\rho_0,u_0)\in H^{k+1}(\mathbb T)\times H^k(\mathbb T)$ with $\min_{\mathbb T}\rho_0>0$. Let $f\in L^2_{\mathrm{loc}}([0,\infty);H^{k-1}(\mathbb T))$. The local construction recalled in \cite[Appendix A]{Ant} provides a unique maximal local strong solution $(\rho,u)$, with maximal existence time $T^*\in(0,\infty]$, such that
\begin{equation}\label{eq:NSK-local-regularity}
\begin{split}
\rho&\in C([0,T^*);H^{k+1}(\mathbb T)),\\
u&\in C([0,T^*);H^k(\mathbb T))
\cap L^2_{\mathrm{loc}}([0,T^*);H^{k+1}(\mathbb T)),
\end{split}
\end{equation}
and $\rho(x,t)>0$ on $\mathbb T\times[0,T^*).$
\\
\\
The additional derivative on $\rho$ reflects the capillary energy and the Korteweg force, which contains third-order density derivatives. In the augmented formulation below, a density gradient and $u$ are estimated at the same Sobolev order.
The friction and external force are lower-order terms and can be included in the same construction. For completeness, we also mention the detailed analysis provided in  \cite[Theorem 2.3]{Pesc} which accounts for both the deterministic and stochastically forced problems. Note that the condition $k>7/2$ used there is due to the stochastic construction, whereas its deterministic counterpart applies for $k>5/2$, as explained in \cite[Remark 2.5(2)]{Pesc}.
\\
\\
The unique local strong solution can be continued as long as $\|\rho^{-1}\|_{L^\infty}$, $\|\rho\|_{H^{k+1}}$ and $\|u\|_{H^k}$ remain bounded. The following theorem shows that the density estimates
of Section \ref{Sec2}, combined with estimates for higher-order
derivatives, imply $T^*=\infty$.
\begin{theorem}\label{thm:NSK-global}
Let $\alpha,\theta\ge0$, $\beta\in\mathbb R$, $\gamma\ge1$,
$\varepsilon>0$, and let $k\ge3$ be an integer.
Assume that
\begin{equation}\label{eq:NSK-global-data}
(\rho_0,u_0)\in H^{k+1}(\mathbb T)\times H^k(\mathbb T),
\qquad
\min_{\mathbb T}\rho_0>0,
\end{equation}
and
\begin{equation*}
f\in L^2_{\mathrm{loc}}
\bigl([0,\infty);H^{k-1}(\mathbb T)\bigr).
\end{equation*}
Suppose that the strong coercivity condition
\begin{equation}\label{scc theorem global wp}
2\alpha-4\le\beta\le2\alpha-1
\end{equation}
holds and that one of the following conditions is satisfied:
\begin{itemize}
\item[(i)]
$0\le\alpha\le\max\{\frac12,\frac{2-\theta}{3}\}$.
\item[(ii)]
$\beta\le\max\{-2,-\alpha-\theta-1\}$.
\end{itemize}
Then system \eqref{NSK damping} admits a unique global
strong solution satisfying, for every $T>0$,
\begin{equation}\label{eq:NSK-global-regularity}
\begin{split}
\rho&\in C([0,T];H^{k+1}(\mathbb T))
       \cap L^2(0,T;H^{k+2}(\mathbb T)),\\
u&\in C([0,T];H^k(\mathbb T))
       \cap L^2(0,T;H^{k+1}(\mathbb T)),
\end{split}
\end{equation}
and
\begin{equation}\label{eq:NSK-global-density}
0<c_T\le\rho(x,t)\le C_T
\qquad\text{on }\mathbb T\times[0,T].
\end{equation}
\end{theorem}
\begin{proof}
Let $(\rho,u)$ be the maximal local strong solution in the class \eqref{eq:NSK-local-regularity}, and fix $T<T^*$. All the estimates below depend only on $T$, the initial data, the coefficients and the norm of $f$ in $L^2(0,T;H^{k-1})$. Proposition \ref{Prop density NSK} gives \eqref{eq:NSK-global-density} under conditions (i) or (ii).
We introduce the following notation:
\begin{equation*}
A=\sqrt{\frac{k(\rho)}{\rho}}\,\partial_x\rho,
\qquad
q(\rho)=\sqrt{\rho k(\rho)},
\qquad
\nu(\rho)=\frac{\mu(\rho)}{\rho}.
\end{equation*}
The above change of variables highlights a skew-symmetric formulation of system \eqref{NSK damping} in the new variables $(A,u)$. In particular, since $\mathcal K=\rho q(\rho)\partial_xA$,
the continuity and velocity equations give
\begin{equation}\label{eq:NSK-global-augmented}
\begin{split}
\partial_tA+\partial_x(q(\rho)\partial_xu)
&=-\partial_x(uA),\\
\partial_tu-\partial_x(\nu(\rho)\partial_xu)
-\partial_x(q(\rho)\partial_xA)
&=-u\partial_xu-\frac{p'(\rho)}{q(\rho)}A
+\frac{\nu(\rho)}{q(\rho)}A\partial_xu\\
&\quad+A\partial_xA-\varepsilon\rho^{\theta-1}u+f.
\end{split}
\end{equation}
By virtue of the skew-symmetric structure, the principal higher-order terms involving $q(\rho)$ cancel in the high-order derivative estimates for $(A,u)$. More precisely, for every $0\le j\le k$, integration by parts gives
\begin{equation}\label{eq:NSK-global-cancellation}
\begin{split}
&-\int_{\mathbb T}\partial_x^jA\,
\partial_x\bigl(q(\rho)\partial_x^{j+1}u\bigr)\,dx+\int_{\mathbb T}\partial_x^ju\,
\partial_x\bigl(q(\rho)\partial_x^{j+1}A\bigr)\,dx=0.
\end{split}
\end{equation}
The remaining terms, including the commutators arising from differentiation of $q(\rho)$, are estimated using interpolation and the viscous dissipation. The energy estimate and the density bounds control $\|\partial_x\rho\|_{L^2}$ uniformly on $[0,T]$. Moreover, the expression for the capillary contribution $\mathcal F_k$ to the BD dissipation gives
\begin{equation*}
\mathcal F_k[\rho]
\ge c_T\|\partial_{xx}\rho\|_{L^2}^2
-C_T\|\partial_x\rho\|_{L^4}^4
\ge \frac{c_T}{2}\|\partial_{xx}\rho\|_{L^2}^2-C_T.
\end{equation*}
The last inequality follows from one-dimensional interpolation and the bound for $\|\partial_x\rho\|_{L^2}$. Thus the BD entropy controls
$\partial_{xx}\rho^{\frac{\alpha+\beta+1}{2}}$ in $L^2(0,T;L^2)$ when $\alpha+\beta+1\ne0$, with $\partial_{xx}\log\rho$ in the case $\alpha+\beta+1=0.$ This second-order density estimate is fundamental to start the higher-order argument, since it gives
$A\in L^2(0,T;H^1)$. Consequently, the energy and BD entropy estimates imply
\begin{equation}\label{eq:NSK-global-base}
A,u\in L^\infty(0,T;L^2)\cap L^2(0,T;H^1).
\end{equation}
This argument also applies at the endpoints of \eqref{scc theorem global wp}. We next estimate higher-order derivatives by defining 
$$Y_j=\|A\|_{H^j}^2+ \|u\|_{H^j}^2.$$
Differentiating \eqref{eq:NSK-global-augmented} once and testing with $(\partial_xA,\partial_xu)$, we use \eqref{eq:NSK-global-cancellation} with $j=1$ to eliminate the principal terms involving $q(\rho).$ Estimating the remaining terms with standard H\"{o}lder and Sobolev arguments we have
\begin{equation*}
\dfrac{d}{dt} Y_1+c_T\|u\|_{H^2}^2
\le C_T\left(1+Y_1+\|f\|_{L^2}^2\right)(1+Y_1).
\end{equation*}
Since $Y_1\in L^1(0,T)$ by \eqref{eq:NSK-global-base}, Gronwall's lemma gives
\begin{equation*}
A,u\in L^\infty(0,T;H^1),
\qquad
u\in L^2(0,T;H^2).
\end{equation*}
The first equation in \eqref{eq:NSK-global-augmented} then yields $\partial_tA\in L^2(0,T;L^2)$.
Indeed,
\begin{equation*}
\partial_tA
=-q(\rho)\partial_{xx}u
-\partial_xq(\rho)\partial_xu
-u\partial_xA-A\partial_xu,
\end{equation*}
and all four terms belong to $L^2(0,T;L^2).$ We next test the velocity equation with $\partial_{xx}A.$ The term involving $\partial_tu$ is handled through
\begin{equation}\label{eq:NSK-global-mixed}
\int_{\mathbb T}\partial_tu\,\partial_{xx}A\,dx
=-\frac{d}{dt}\int_{\mathbb T}\partial_xu\,\partial_xA\,dx
-\int_{\mathbb T}\partial_{xx}u\,\partial_tA\,dx.
\end{equation}
Since $q(\rho)$ is bounded away from zero, integration in time and Young's inequality give
\begin{equation*}
A\in L^2(0,T;H^2).
\end{equation*}
Thus $Y_2\in L^1(0,T)$. We now differentiate \eqref{eq:NSK-global-augmented} twice and test with $(\partial_{xx}A,\partial_{xx}u)$.
The principal terms involving $q(\rho)$ cancel by
\eqref{eq:NSK-global-cancellation} with $j=2$, and we obtain
\begin{equation*}
\dfrac{d}{dt} Y_2+c_T\|u\|_{H^3}^2
\le C_T\left(1+Y_2+\|f\|_{H^1}^2\right)(1+Y_2).
\end{equation*}
Another application of Gronwall's lemma yields
\begin{equation}\label{eq:NSK-global-H2}
A,u\in L^\infty(0,T;H^2),
\qquad
u\in L^2(0,T;H^3).
\end{equation}
These bounds close the estimates at every higher order. For each $3\le j\le k$, we apply $\partial_x^j$ to \eqref{eq:NSK-global-augmented} and test with $(\partial_x^jA,\partial_x^ju)$.
The same cancellation \eqref{eq:NSK-global-cancellation} eliminates the principal terms involving $q(\rho)$ at order $j$. After estimating the commutators and the remaining terms, we find
\begin{equation*}
\dfrac{d}{dt}Y_j+c_T\|u\|_{H^{j+1}}^2
\le C_T\left(1+\|u\|_{H^3}^2\right)(1+Y_j)
+C_T\|f\|_{H^{j-1}}^2.
\end{equation*}
Here the constants also depend on the bounds in
\eqref{eq:NSK-global-H2}. In particular, the coefficients involving
$\partial_xA$ and $\partial_xu$ are bounded, while $\|\partial_{xx}u\|_{L^\infty}$ is integrable in time because of the embedding $H^3(\mathbb T) \hookrightarrow W^{2,\infty}(\mathbb T)$. Gronwall's lemma therefore gives
\begin{equation*}
A,u\in L^\infty(0,T;H^k),
\qquad
u\in L^2(0,T;H^{k+1}).
\end{equation*}
To obtain the $L^2(0,T;H^{k+2})$ bound for $\rho$, we use the same mixed estimate at order $k$. The equation for $A$ in \eqref{eq:NSK-global-augmented} and the preceding bounds give $\partial_tA\in L^2(0,T;H^{k-1})$. Applying $\partial_x^{k-1}$ to the velocity equation and testing with $\partial_x^{k+1}A$, we integrate by parts as in \eqref{eq:NSK-global-mixed} to obtain
\begin{equation*}
\int_{\mathbb T}\partial_x^{k-1}\partial_tu\,\partial_x^{k+1}A\,dx
=-\frac{d}{dt}\int_{\mathbb T}\partial_x^ku\,\partial_x^kA\,dx
-\int_{\mathbb T}\partial_x^{k+1}u\,\partial_x^{k-1}\partial_tA\,dx.
\end{equation*}
After integration in time, the first term is controlled by the $L^\infty(0,T;H^k)$ bounds for $u$ and $A$, and the second by the $L^2$ bounds for $\partial_x^{k+1}u$ and $\partial_x^{k-1}\partial_tA$. Estimating the commutators and the forcing as in \eqref{eq:NSK-global-mixed}, we obtain $A\in L^2(0,T;H^{k+1})$. The definition of $A$ and the upper and lower bounds on $\rho$ then give
\begin{equation*}
\rho\in L^\infty(0,T;H^{k+1})\cap L^2(0,T;H^{k+2}).
\end{equation*}
The equations also give $\partial_tu\in L^2(0,T;H^{k-1})$ and $\partial_t\rho\in L^2(0,T;H^k)$. Together with the spatial estimates, these bounds yield the continuity in time stated in \eqref{eq:NSK-global-regularity}. If $T^*$ were finite, the preceding estimates would bound $\|\rho(t)^{-1}\|_{L^\infty}$, $\|\rho(t)\|_{H^{k+1}}$ and $\|u(t)\|_{H^k}$ uniformly for $0\le t<T^*$. The local continuation criterion would then extend the solution beyond $T^*$, contradicting its maximality. 
\end{proof}
\noindent
\textbf{Declaration of AI use.} The author used GPT-6.2 Astra to correct typographical and grammatical errors, check the presentation and double-check the proofs.  In addition, GPT-6.2 Astra was used to improve an earlier version of Proposition \ref{active potential extended thm}.
The mathematical ideas and strategy of proofs belong to the author, who takes full responsibility for the content of this paper.
\par\smallskip\noindent
\textbf{Data Availability.} The author declares that data sharing is not applicable to this article as no datasets were generated or analysed.
\par\smallskip\noindent
\textbf{Conflict of interest.} The author declares that there is no conflict of interest.

\end{document}